\documentclass[leqno,12pt]{article} %twoside is a good option.

\usepackage{amsmath}
\usepackage{amssymb}
\usepackage{pst-all}

\newtheorem{theorem}{Theorem}[section]

\newtheorem{definition}[theorem]{Definition}

\newtheorem{note}[theorem]{Note}

\newtheorem{proposition}[theorem]{Proposition}
\newtheorem{corollary}[theorem]{Corollary}

\newtheorem{lemma}[theorem]{Lemma}

\newcommand{\R}{\mathbb{R}}
\newcommand{\N}{\mathbb{N}}

\newcommand{\qed}{\ensuremath{\quad \blacksquare} \vspace{\baselineskip}}

\newcommand{\Diag}{\mathrm{Diag\,}}
\newcommand{\diag}{\mathrm{diag\,}}

\newcommand{\signum}{\mathrm{sign\,}}

\begin{document}
\allowdisplaybreaks
%***********************************************************************
\title{\textbf{Interactions Between the Generalized Hadamard
Product and the Eigenvalues of Symmetric Matrices}}
\author{\textbf{
Hristo S.\ Sendov\thanks{Department of Mathematics and Statistics,
University of Guelph, Guelph, Ontario, Canada N1G 2W1. Email:
\texttt{hssendov\char64 uoguelph.ca}. Research supported by
NSERC.}
 }}

\maketitle

\begin{abstract}
In \cite{Sendov2003a} a notion of generalized Hadamard product was
introduced. We show that when certain kinds of tensors interact
with the eigenvalues of symmetric matrices the resulting formulae
can be nicely expressed using the generalized Hadamard product and
two simple linear operations on the tensors. The Calculus-type
rules developed here will be used in \cite{Sendov2003c} to
routinize, to a large extend, the differentiation of spectral
functions. We include all the necessary definitions and results
from \cite{Sendov2003a} on the generalized Hadamard product to
make the reading self-contained.
\end{abstract}

\noindent {\bf Keywords:}  spectral functions, eigenvalues,
Hadamard product, tensor analysis, eigenvalue optimization,
symmetric function, perturbation theory.

\noindent {\bf Mathematics Subject Classification (2000):}
primary: 49R50; 47A75, secondary: 15A18; 15A69.

%**********************************************************************
\section{Introduction}

The aim of this paper is to develop some analytic tools that, we
believe, will fully describe the formula for the higher order
derivatives of {\it spectral functions} in terms of the underlying
{\it symmetric function}. We say that a real-valued function $F$,
on a symmetric matrix argument, is spectral if it has the
following invariance property:
$$
F(UXU^T)=F(X),
$$
for every symmetric matrix $X$ in its domain and every orthogonal
matrix $U$. When $U$ varies freely over the orthogonal matrices
the invariants of the product $UXU^T$ are the eigenvalues of the
matrix $X$. Therefore, de facto, the function $F$ depends only on
the set of eigenvalues of $X$. The restriction of $F$ to the
subspace of diagonal matrices defines (almost) a function
$f(x):=F(\Diag x)$ on a vector argument $x \in \R^n$. It is easy
to see that $f: \R^n \rightarrow \R$ has the property
$$
f(x)=f(Px) \mbox{  for any permutation matrix  }  P \mbox{ and any
} x \in \mbox{domain} \, f.
$$
We call such functions {\it symmetric}. It is not difficult to see
that $F(X) = f(\lambda(X))$, where $\lambda(X)$ is the vector of
eigenvalues of $X$.

One of the main questions in the theory of spectral functions is
what smoothness properties of the symmetric function $f$ are
inherited by $F$. The difficulties arise from the fact that the
eigenvalue map, $\lambda(X)$, don't depend smoothly on its
argument $X$. Even in domains where they are smooth, it is
difficult to organize the differentiation process so that the
result is in as closed form as possible.

One of the fist results in this direction (see \cite{Lewis:1994b})
showed that $F$ is (continuously) differentiable at a matrix $A$
if and only if $f$ is, at the vector $\lambda(A)$. The formula for
the gradient is compact and easy to understand. Next,
\cite{LewisSendov:2000a}, showed that $F$ is twice (continuously)
differentiable if, and only if, $f$ is at $\lambda(A)$. The
consideration of variety of different subcases made the
differentiation laborious and the formula for the Hessian of $F$
takes some effort to get comfortable with. Following the
developments in \cite{LewisSendov:2000a}, one can see that an
attempt to compute the third, or higher, derivatives of $F$, will
result in a number of subcases that will quickly become
unmanageable. That is why, deriving a formula for the higher
derivatives (and in the process proving that they exist) requires
a language that handles all the cases in a structured way, and
allows easy to work with calculus rules. In \cite{Sendov2003a} we
proposed such a language based on the idea of generalizing the
Hadamard product between two matrices to a tensor-valued product
between $k$ matrices, $k \ge 1$. This paper is a continuation of
the work there. While, \cite{Sendov2003a} emphasizes on
multi-linear algebra and combinatorial aspects of the generalized
Hadamard product, our current work deals mainly with the calculus
of the generalized Hadamard product that is related to the
eigenvalues of symmetric matrices.

It is likely that high-powered analytical methods will be able to
show directly that $F$ is $k$ times (continuously) differentiable
if and only if $f$ is, see \cite{GrabovskyHijabRivin:2003}. Our
approach aims to give a constructive procedure how to, knowing the
$k$-th derivative of $F$, practically compute the $(k+1)$-st. The
precise description of the directional expansion of the
eigenvalues of a symmetric matrix, when the entries of the matrix
depend only on one scalar parameter (see \cite{kato:1976full},
\cite{kemble:1958}) finds far reaching practical applications in
areas ranging from optimization to quantum mechanics. Formulae for
the derivatives of spectral functions will naturally include, as a
special case, the directional derivatives, when the symmetric
matrix depends on one argument.

The paper is organized as follows. In the next section we
introduce the necessary notation and definitions. The background
definitions and results from \cite{Sendov2003a} that will be
needed are given in Section~\ref{sec-2}. All this aims to make the
reading as self-contained as possible. In Section~\ref{refinment}
we reexamine Lemma~2.4 from \cite{LewisSendov:2000a}. We distill
the essential eight parts of the statement of the lemma down to
two equations: the first is nothing more than a spectral
decomposition of a block-diagonal symmetric matrix, while the
second is a (strong) first-order expansion for dot products
between (parts of) eigenvectors. The section after that contains
the calculus results. The main theorems of that section are
Theorem~\ref{dec14b} and Theorem~\ref{dec15a}. They investigate
how the interaction between {\it block-constant tensors} (see
below) and eigenvalues affect the generalized Hadamard product. We
use these tools in \cite{Sendov2003c} to derive a computable
formula for the higher-order derivatives of any spectral function
at a symmetric matrix $X$ with distinct eigenvalues, as well as
for the derivatives of separable spectral functions at an
arbitrary symmetric matrix. (Separable spectral functions are
those arising from symmetric functions $f(x) = g(x_1) + \cdots +
g(x_n)$ for some function $g$ on a scalar argument.)

%**********************************************************************
\section{Notation and definitions}
\label{NAPR}

In what follows, $M^n$ will denote the Euclidean space of all $ n
\times n$ real matrices with inner product $\langle X,Y \rangle =
\mbox{tr}\,(XY^T)$, and the subspace of $n \times n$ symmetric
matrices will be denoted by $S^n$. For $A \in S^n$,
$\lambda(A)=(\lambda_1(A), ...,\lambda_n(A))$ will be the vector
of its eigenvalues ordered in nonincreasing order.  By $O^n$ and
$P^n$ we will denote the set of all $n \times n$ orthogonal and
permutation matrices respectively. By $\N_k$ we will denote the
set $\{1,2,...,k\}$. For any vector $x$ in $\R^n$, Diag$\,x$ will
denote the diagonal matrix with the vector $x$ on the main
diagonal, and diag$\,: M^n \rightarrow \R^n$ will denote its
conjugate operator, defined by diag$\,(X)=(x_{11},...,x_{nn})$. By
$\R_{\downarrow}^n$ we denote the cone of all vectors $x$ in
$\R^n$ such that $x_1 \ge x_2 \ge \cdots \ge x_n$. In the whole
paper $\{M_m \}_{m=1}^{\infty}$ will denote a sequence of
symmetric matrices converging to $0$, $\{U_m \}_{m=1}^{\infty}$
will denote a sequence of orthogonal matrices. We describe sets in
$\R^n$ and functions on $\R^n$ as {\em symmetric} if they are
invariant under coordinate permutations. We denote the gradient of
$f$ by $\nabla f$, and the Hessian by $\nabla^2 f$. In general, if
the function $f$ is $k$ times differentiable, then $\nabla^k
f(\mu)$ will denote its $k$-th order differential at the point
$\mu$. It can be viewed as a $k$-dimensional tensor on $\R^n$.

\begin{definition} \rm
A {\it $k$-tensor, $T$, on $\R^n$} is a map from $\R^n \times
\cdots \times \R^n$ ($k$-times) to $\R$ that is linear in each
argument separately. The set of all $k$-tensors on $\R^n$ will be
denoted by $T^{k,n}$. The value of the $k$-tensor, $T$, at
$(h_1,...,h_k)$ will be denoted by $T[h_1,...,h_k]$. The tensor is
called {\it symmetric} if for any permutation, $\sigma$, on $\N_k$
it satisfies
$$
T[h_{\sigma(1)},...,h_{\sigma(k)}] = T[h_1,...,h_k],
$$
for any $h_1,...,h_k \in \R^n$.
\end{definition}

We denote the standard basis in $\R^n$ by $e^1,e^2,...,e^n$. For
an arbitrary $k$-tensor, $T$, and any $k$-tuple of integers from
$\N_n$, $(i_1,...,i_k)$, we denote its $(i_1,...,i_k)$-th element
by $T^{i_1...i_k}$. (Matrices will be viewed as $2$-tensors and
vectors as $1$-tensors.) If $T \in T^{k,n}$ and $h \in \R^n$, then
for brevity throughout the paper, we denote by $T[h]$ the
$(k-1)$-tensor on $\R^n$ given by $T[\cdot,...,\cdot,h]$.
Similarly for $T[M]$, when $T$ is a $k$-tensor on $M^n$ and $M\in
M^n$.

For a permutation matrix $P \in P^n$ we say that $\sigma : \N_n
\rightarrow \N_n$ is its corresponding permutation map and write
$P \leftrightarrow \sigma$ if for any $h \in \R^n$ we have
$Ph=(h_{\sigma(1)},...,h_{\sigma(n)})^T$ or, in other words,
$P^Te^i = e^{\sigma(i)}$ for all $i=1,...,n$. The symbol
$\delta_{ij}$ will denote the Kroneker delta. It is equal to one
if $i=j$ and zero otherwise.

Any vector $\mu \in \R^n$ defines a partition of $\N_n$ into
disjoint {\it blocks}, where integers $i$ and $j$ are in the same
block if, and only if, $\mu_i=\mu_j$. Whenever $\mu$ is a vector
in $\R^n_{\downarrow}$ we make the convention that
$$
\mu_1= \cdots = \mu_{k_1}> \mu_{k_1+1}= \cdots = \mu_{k_2} >
\mu_{k_2+1} \cdots \mu_{k_r}, \hspace{1.5cm} (k_0=0, \, k_r=n),
$$
with corresponding partition
\begin{equation}
\label{jen} I_1:=\{1,2,...,\iota_1
\},~I_2:=\{\iota_1+1,\iota_1+2,...,\iota_2
\},...,~I_r:=\{\iota_{r-1}+1,...,\iota_r \}.
\end{equation}

For arbitrary vector $\mu$ the blocks it determines need not
contain consecutive integers. Thus, we agree that the block
containing the integer $1$ will be the first block, $I_1$, the
block containing the smallest integer that is not in $I_1$ will be
the second block, $I_2$, and so on. This naturally enumerates all
 the blocks, and in general, $\iota_l$ will denote the largest
 integer in $I_l$ for all $l = 1,...,r$. Also, $r$ will denote the number of
blocks determined by $\mu$. For any two integers, $i,j \in \N_n$
we will say that they are {\it equivalent (with respect to $\mu$)}
and write $i \sim j$ (or $i \sim_{\mu} j$) if $\mu_i = \mu_j$,
that is, if they are in the same block. Two $k$-indexes
$(i_1,...,i_k)$ and $(j_1,...,j_k)$ are called {\it equivalent} if
$i_l \sim j_l$ for all $l=1,2,...,n$, and we will write
$$
(i_1,...,i_k) \sim (j_1,...,j_k).
$$

\begin{definition} \rm
Given a vector $\mu \in \R^n$, we say that a $k$-tensor, $T$, is
{\it ($\mu$-)block-constant} if $T^{i_1...i_k}=T^{j_1...j_k}$
whenever $(i_1,...,i_k) \sim (j_1,...,j_k)$.

A $k$-tensor valued map, $\mu \in \R^n \rightarrow
\mathcal{F}(\mu) \in T^{k,n}$, is {\it block-constant} if
$\mathcal{F}(\mu)$ is $\mu$-block-constant for every $\mu$.
\end{definition}

The following elementary lemma motivates the definitions following
it. It follows from applying the chain rule to the equality
$f(\mu)=f(P\mu)$.

\begin{lemma}
\label{Block:l} Let $f:\R^n \rightarrow \R$ be a symmetric
function, $k$ times differentiable at the point $\mu \in \R^n$,
and let $P$ be a permutation matrix such that $P\mu = \mu$.  Then
\begin{enumerate}
\item \label{one:block:l} $\nabla f(\mu) = P^T \nabla f(\mu)$,
\item \label{two:block:l} $\nabla^2 f(\mu) = P^T \nabla^2 f(\mu)
P$, and in general \item \label{three:block:l}$\nabla^s
f(\mu)[h_1,...,h_s] = \nabla^s f(\mu)[Ph_1,...,Ph_s]$, for any
$h_1,...,h_s \in \R^n$, and $s \in \N_k$.
\end{enumerate}
\end{lemma}

\begin{definition} \rm
\label{mu-symm:defn}  Given a vector $\mu \in \R^n$, we will say
that $T \in T^{k,n}$ is {\it $\mu$-symmetric} if for any
permutation $P \in P^n$, such that $P\mu = \mu$, we have
$$
T[Ph_{1},...,Ph_{k}] = T[h_1,...,h_k], \mbox{ for any }
h_1,...,h_k \in \R^n.
$$
A $k$-tensor valued map, $\mu \in \R^n \rightarrow
\mathcal{F}(\mu) \in T^{k,n}$, is {\it $\mu$-symmetric} if for
every $\mu \in \R^n$ and permutation matrix $P$ we have
$$
\mathcal{F}(P\mu)[Ph_1,...,Ph_k] =
\mathcal{F}(\mu)[h_{1},...,h_{k}], \mbox{ for any } h_1,...,h_k
\in \R^n.
$$
\end{definition}

Clearly, every $\mu$-block-constant tensor is $\mu$-symmetric, the
opposite is not true. There is a slight abuse of terminology since
$\mu$-symmetric means different things for a $k$-tensor and for a
$k$-tensor valued map. If the map $\mu \in \R^n \rightarrow
\mathcal{F}(\mu) \in T^{k,n}$ is $\mu$-symmetric, then for a fixed
$\mu$ the tensor $\mathcal{F}(\mu)$ is $\mu$-symmetric. This makes
sure there will be no confusion.

By Lemma~\ref{Block:l}, for any differentiable enough, symmetric
function $f : \R^n \rightarrow \R$ the mapping $\mu \in \R^n
\rightarrow \nabla f(\mu) \in \R^n$ is a $\mu$-symmetric,
$\mu$-block-constant, $1$-tensor valued mapping. In general, for
every $s \in \N_k$ the mapping $\mu \in \R^n \rightarrow \nabla^s
f(\mu)$ is a $\mu$-symmetric, $s$-tensor-valued map, and if
continuous, then every tensor is also symmetric.

We conclude this section with the following easy lemma.

\begin{lemma}
\label{sec1-lastlem} If a $k$-tensor valued map, $\mu \in \R^n
\rightarrow T(\mu) \in T^{k,n}$, is $\mu$-symmetric and
differentiable, then its differential is also $\mu$-symmetric.
\end{lemma}

\begin{proof}
We use the first-order Taylor expansion formula. Let $v_m$ be a
sequence of vectors approaching zero such that $v_m/\|v_m\|$
approaches $h$ as $m \rightarrow \infty$.
\begin{align*}
T(\mu+v_m)[h_1,...,h_k] = T(\mu)[h_1,...,h_k] + \nabla
T(\mu)[h_1,...,h_k, v_m] + o(\|v_m\|).
\end{align*}
On the other hand, for any permutation $P$ we have
\begin{align*}
T(\mu+v_m)[h_1,...,h_k] &= T(P\mu+Pv_m)[Ph_1,...,Ph_k] \\
&= T(P\mu)[Ph_1,...,Ph_k] + \nabla
T(P\mu)[Ph_1,...,Ph_k, Pv_m] + o(\|Pv_m\|) \\
&= T(\mu)[h_1,...,h_k] + \nabla T(P\mu)[Ph_1,...,Ph_k, Pv_m] +
o(\|v_m\|).
\end{align*}
Subtracting the two equalities, dividing by $\|v_m\|$ and letting
$m$ go to infinity, we get
$$
\nabla T(P\mu)[Ph_1,...,Ph_k, Ph] = T(\mu)[h_1,...,h_k, h].
$$
Since the vectors $h_1$,...,$h_k$, and $h$ are arbitrary, the
result follows. \hfill \qed
\end{proof}

\section{Generalized Hadamard product}
\label{sec-2}

In this section we will quote briefly several definitions and
results from \cite{Sendov2003a} that are crucial for the
development in this work. Recall that the Hadamard product of two
matrices $A=[A^{ij}]$ and $B=[B^{ij}]$ of the same size is the
matrix of their element-wise product $A \circ B = [A^{ij}B^{ij}]$.
The standard basis on the space $M^n$ is given by the set
$\{H_{pq} \in M^n \, | \, H_{pq}^{ij} = \delta_{ip}\delta_{jq}
\mbox{ for all } i,j \in \N_n \}$, where $\delta_{ij}$ is the
Kronecker delta function, equal to one if $i=j$, and zero
otherwise.

\begin{definition} \rm
\label{tripleHad} For each permutation $\sigma$ on $\N_k$, we
define {\it $\sigma$-Hadamard} product between $k$ matrices to be
a $k$-tensor on $\R^n$ as follows. Given any $k$ basic matrices
$H_{p_1q_1}$, $H_{p_2q_2}$,...,$H_{p_kq_k}$
\begin{align*}
 (H_{p_1q_1} \circ_{\sigma} H_{p_2q_2}
\circ_{\sigma} \cdots \circ_{\sigma} H_{p_kq_k})^{i_1i_2...i_k} &=
\left\{
\begin{array}{ll}
1, \hspace{-0.2cm} & \mbox{ if } i_s=p_s=q_{\sigma(s)}, \forall s=1,...,k, \\
0, \hspace{-0.2cm} & \mbox{ otherwise. }
\end{array}
\right.
\end{align*}
\end{definition}

Extend this product to a multi-linear map on $k$ matrix arguments:
\begin{align}
\label{ext-by-lin} (H_{1}\circ_{\sigma} H_{2} \circ_{\sigma}
\cdots \circ_{\sigma} H_k)^{i_1i_2...i_k} &=
H_1^{i_1i_{\sigma^{-1}(1)}}\cdots H_k^{i_ki_{\sigma^{-1}(k)}}.
\end{align}

Let $T$ be an arbitrary $k$-tensor on $\R^n$ and let $\sigma$ be a
permutation on $\N_k$. We define $\Diag^{\sigma} T$ to be a
$2k$-tensor on $\R^n$ in the following way
\begin{align*}
(\Diag^{\sigma} T)^{\substack{i_1...i_k \\ j_1...j_k}} = \left\{
\begin{array}{ll}
T^{i_1...i_k}, & \mbox{ if } i_s = j_{\sigma(s)}, \forall s = 1,...,k, \\
0, & \mbox{ otherwise. }
\end{array}
\right.
\end{align*}
Notice that any $2k$-tensor, $T$, on $\R^n$ can naturally be
viewed as a $k$-tensor on $M^n$ in the following way
$$
T[H_1,...,H_k]= \sum_{p_1,q_1=1}^n \cdots \sum_{p_k,q_k=1}^n
T^{\substack{p_1...p_k \\ q_1...q_k}} H_1^{p_1q_1} \cdots
H_k^{p_kq_k}.
$$
Define dot product between two tensors in $T^{k,n}$ in the usual
way:
$$
\langle T_1, T_2\rangle = \sum_{p_1,...,p_k = 1}^n
T_1^{p_1...p_k}T_2^{p_1...p_k}.
$$
 We define an action (called {\it conjugation}) of the orthogonal group $O^n$ on the
space of all $k$-tensors on $\R^n$. For any $k$-tensor, $T$, and
$U \in O^n$ this action will be denoted by $UTU^T \in T^{k,n}$:
\begin{equation}
\label{kur0} (UTU^T)^{i_1...i_k} = \sum_{p_1 = 1 }^{n} \cdots
\sum_{p_k =1}^{n} \Big( T^{p_1...p_k} U^{i_{1}p_{1}} \cdots
U^{i_{k}p_{k}} \Big).
\end{equation}
It was shown in \cite{Sendov2003a} that this action is norm
preserving and associative: $V(UTU^T)V^T = (VU)T(VU)^T$ for all
$U,V \in O^n$.

The $\Diag^{\sigma}$ operator, the $\sigma$-Hadamard product, and
conjugation by an orthogonal matrix are connected by the following
formula, see \cite{Sendov2003a}.
\begin{theorem}
\label{jen3a} For any $k$-tensor $T$, any matrices
$H_1$,...,$H_k$, any orthogonal matrix $V$, and any permutation
$\sigma$ in $P^k$ we have the identity
\begin{equation}
\label{gen:1} \langle T, \tilde{H}_1 \circ_{\sigma} \cdots
\circ_{\sigma} \tilde{H}_k \rangle = \big(V(\Diag^{\sigma}
T)V^T\big)[H_1,...,H_k],
\end{equation}
where $\tilde{H}_i = V^TH_iV$, $i=1,...,k$.
\end{theorem}

\begin{lemma}
\label{dot-prod} Let $T$ be a $k$-tensor on $\R^n$, and $H$ be a
matrix in $M^n$. Let $H_{i_1j_1}$,...,$H_{i_{k-1}j_{k-1}}$ be
basic matrices in $M^n$, and let $\sigma$ be a permutation on
$\N_k$. Then the following identities hold.
\begin{enumerate}
\item If $\sigma^{-1}(k)=k$, then
\begin{align*}
\langle T, H_{i_1j_1} &\circ_{\sigma} \cdots \circ_{\sigma}
H_{i_{k-1}j_{k-1}} \circ_{\sigma} H \rangle =
\Big(\prod_{t=1}^{k-1} \delta_{i_tj_{\sigma(t)}} \Big)
\sum_{t=1}^n T^{i_1 \dots i_{k-1}t}H^{tt}.
\end{align*}
\item If $\sigma^{-1}(k) = l$, where $l \not= k$, then
\begin{align*}
\langle T, H_{i_1j_1} &\circ_{\sigma} \cdots \circ_{\sigma}
H_{i_{k-1}j_{k-1}} \circ_{\sigma} H \rangle =
\Big(\prod_{\substack{t=1 \\ t \not= l}}^{k-1}
\delta_{i_tj_{\sigma(t)}} \Big) T^{i_1 \dots
i_{k-1}j_{\sigma(k)}}H^{j_{\sigma(k)}i_{\sigma^{-1}(k)}}.
\end{align*}
\end{enumerate}
\end{lemma}

\section{A refinement of a perturbation result for eigenvectors}
\label{refinment}

The main limiting tool in \cite{LewisSendov:2000a} was Lemma~2.4.
The statement of the lemma was broken down into nine different
parts, and that lead to the consideration of variety of cases when
deriving the formula for the Hessian of spectral functions. For
general situations, that we aim to tackle later, such case studies
will quickly become unmanageable. That is why the goal of this
section is to transform Lemma~2.4 from \cite{LewisSendov:2000a}
into a form more suitable for computations. We begin with a lemma
(the proof is a simple combination of Lemma~5.10 in
\cite{Lewis:1996} and Theorem~3.12 in \cite{HUYe:1995}) that will
alow us to define some of the necessary notation.

% We begin this section by quoting a lemma, .

\begin{lemma}
\label{lemma:1} Let $\mu \in \R^n_{\downarrow}$ and let
(\ref{jen}) be the partition defined by $\mu$. For any sequence of
symmetric matrices $M_m \rightarrow 0$ we have that
\begin{equation}
\label{dec15a:eqn} \lambda(\mbox{\rm Diag}\, \mu + M_m)^T = \mu^T
+ \left( \lambda(X_1^TM_m X_1)^T,...,\lambda(X_r^TM_mX_r)^T
\right)^T + o(\|M_m \|),
\end{equation}
where $X_l := [e^{k_{l-1}+1},...,e^{k_l}], \, \mbox{ for all } \,
l=1,...,r$.
\end{lemma}

 We need some additional notation that
will be used later as well. For a fixed $\mu \in \R^n$ and any
square matrix $H$, we define
\begin{align*}
H^{ij}_{\mbox{\scriptsize \rm in}} &= \left\{ \begin{array}{ll}
H^{ij}, &
\mbox{ if } i \sim j, \\
0, & \mbox{ otherwise, }
\end{array} \right. \\
H^{ij}_{\mbox{\scriptsize \rm out}} &= \left\{ \begin{array}{ll}
H^{ij}, &
\mbox{ if } i \not\sim j,  \\
0, & \mbox{ otherwise. }
\end{array} \right.
\end{align*}
In other words, $H_{\mbox{\scriptsize \rm in}}$ extracts the
diagonal blocks from the matrix $H$ and puts zeros everywhere
else, while $H_{\mbox{\scriptsize \rm out}}$ extracts the off
diagonal blocks and fills the diagonal blocks with zeros. Clearly
for any matrix $H$ we have
$$
H = H_{\mbox{\scriptsize \rm in}} + H_{\mbox{\scriptsize \rm
out}}.
$$

Throughout the whole paper, we denote
\begin{equation}
\label{dec15b:eqn} h_m := \left( \lambda(X_1^TM_m
X_1)^T,...,\lambda(X_r^TM_mX_r)^T \right)^T.
\end{equation}
If also $M_m/\|M_m\|$ converges to $M$ as $m$ goes to infinity,
since the eigenvalues are continuous functions, we can define
\begin{equation}
\label{limh} h := \lim_{m \rightarrow \infty} \frac{h_m}{\|M_m\|}
= \left( \lambda(X_1^TM X_1)^T,...,\lambda(X_r^TMX_r)^T \right)^T.
\end{equation}
We reserve the symbols $h_m$ and $h$ to denote the above two
vectors throughout the paper, unless stated otherwise. With this
notation Lemma~\ref{lemma:1} says that if $M_m \rightarrow 0$,
then
\begin{equation}
\label{kamen} \lambda(\Diag \mu + M_m)^T = \mu^T + h_m +
o(\|M_m\|).
\end{equation}

Below is the main result of this section.

\begin{theorem}
\label{important} Let $\{M_m \}$ be a sequence of symmetric
matrices converging to 0, such that $M_m/\|M_m \|$ converges to
$M$. Let $\mu$ be in $\R^n_{\downarrow}$ and $U_m \rightarrow U
\in O^n$ be a sequence of orthogonal matrices such that
\begin{equation*}
\mbox{\rm Diag}\, \mu + M_m = U_m \big( \mbox{\rm Diag}\,
\lambda(\mbox{\rm Diag}\, \mu + M_m) \big) U_m^T, \, \, \mbox{ for
all } \, \, m=1,2,....
\end{equation*}
Then:
\begin{enumerate}
\item \label{boza1} The orthogonal matrix $U$ has the form
$$U= \left( \begin{array}{cccc}
V_1 & 0   & \cdots & 0 \\
0   & V_2 & \cdots & 0 \\
\vdots & \vdots & \ddots & \vdots \\
0 & 0 & \cdots & V_r
\end{array} \right),
$$
where $V_l$ is an orthogonal matrix with dimensions $|I_l| \times
|I_l|$ for all $l$.
\item \label{boza2} The following identity holds
\begin{equation}
\label{diag} U^T M_{\mbox{\scriptsize \rm in}} U = \Diag h,
\end{equation}
\item \label{boza3} For any indexes $i \in I_l$, $j \in I_s$, and
$t \in \{1,...,r\}$ we have the (strong) first-order expansion
\begin{equation}
\label{exp} \sum_{p \in I_t} U_m^{ip}U_m^{jp} = \delta_{ij}
\delta_{lt} + \frac{\delta_{lt} - \delta_{st}}{\mu_{i} - \mu_{j}}
M^{ij} \|M_m\| +o(\|M_m\|),
\end{equation}
with the understanding that the fraction is zero whenever
$\delta_{lt} = \delta_{st}$ no matter what the denominator is.
\end{enumerate}
\end{theorem}

\begin{proof}
This lemma is essentially Lemma~2.4 in \cite{LewisSendov:2000a}.
Indeed, Part~\ref{boza1} is \cite[Lemma~2.4
Part~(i)]{LewisSendov:2000a}, and Part~\ref{boza2} is an aggregate
version of Parts~(iv) and (vii) from there as well. To prove
Part~\ref{boza3} we consider several cases.
\begin{description}
\item[Case 1.] If $i=j \in I_l$ and $t=l$, then
Formula~(\ref{exp}) becomes $\displaystyle\sum_{p \in I_l}
\big(U_m^{ip}\big)^2 = 1 + o(\|M_m\|),$ which is exactly Part~(ii)
of Lemma~2.4 in \cite{LewisSendov:2000a}. \item[Case 2.] If $i
\not= j \in I_l$ and $t=l$, then Formula~(\ref{exp}) becomes
$\displaystyle\sum_{p \in I_l} \big(U_m^{ip}\big)^2 = o(\|M_m\|),$
which is exactly Part~(vi) of Lemma~2.4 in
\cite{LewisSendov:2000a}. \item[Case 3.] If $i \not= j \in I_l$
and $t \not= l$, then Formula~(\ref{exp}) becomes
$\displaystyle\sum_{p \in I_t} U_m^{ip}U_m^{jp} = o(\|M_m\|),$
which is a consequence of Part~(v) of Lemma~2.4 in
\cite{LewisSendov:2000a}. \item[Case 4.] If $i \in I_l$, $j \in
I_s$, with $l \not= s \not= t \not= l$, then Formula~(\ref{exp})
becomes $\displaystyle\sum_{p \in I_t} U_m^{ip}U_m^{jp} =
o(\|M_m\|),$ which is a consequence of Part~(viii) of Lemma~2.4 in
\cite{LewisSendov:2000a}. \item[Case 5.] If $i \in I_l$, $j \in
I_s$, with $l \not= s$ and $t=l$, then Formula~(\ref{exp}) becomes
$$
\sum_{p \in I_t} U_m^{ip}U_m^{jp} =  \frac{1}{\mu_{i} - \mu_{j}}
M^{ij} \|M_m\| +o(\|M_m\|).
$$
This formula requires a proof. It will be presented together with
the proof of the next, last, case below. \item[Case 6.] If $i \in
I_l$, $j \in I_s$, with $l \not= s$ and $t=s$, then
Formula~(\ref{exp}) becomes
$$
\sum_{p \in I_t} U_m^{ip}U_m^{jp} =  -\frac{1}{\mu_{i} - \mu_{j}}
M^{ij} \|M_m\| +o(\|M_m\|).
$$
We now show that the expressions in both Case~5 and Case~6 are
valid. Recall that Part~(ix) from Lemma~2.4 in
\cite{LewisSendov:2000a} says that in case when $i \in I_l$, $j
\in I_s$ with $l \not= s$, we have
$$
\lim_{m \rightarrow \infty} \bigg( \mu_{k_l} \frac{\sum_{p \in
I_l}U_m^{ip}U_m^{jp} }{\|M_m\|} + \mu_{k_s} \frac{\sum_{p \in
I_s}U_m^{ip}U_m^{jp} }{\|M_m\|} \bigg) = M^{ij}.
$$
Introduce the notation
$$
\beta^l_m := \frac{\sum_{p \in I_l}U_m^{ip}U_m^{jp}}{\|M_m\|},
\hspace{1cm} \mbox{ for all } \, \, \, l=1,2,...,r,
$$
and notice that
$$
\sum_{l=1}^r \beta_m^l = 0, ~~\mbox{ for all } \, m,
$$
because $U_m$ is an orthogonal matrix and the numerator of the
above sum is the product of its $i$-th and the $j$-th row. Next,
Case~4 above says that
$$
\lim_{m \rightarrow \infty} \sum_{t \not= l,s} \beta_m^t = 0,
$$
so
$$
\lim_{m \rightarrow \infty} (\beta_m^l+\beta_m^s) = 0.
$$
For arbitrary reals $a$ and $b$ we compute
\begin{align*}
( a \beta_m^l + b \beta_m^s ) - \frac{ a - b}{\mu_{k_l}-\mu_{k_s}}
\left(\mu_{k_l} \beta_m^l + \mu_{k_s} \beta_m^s \right) =
(\beta_m^l+\beta_m^s) \frac{ b \mu_{k_l} - a
\mu_{k_s}}{\mu_{k_l}-\mu_{k_s}} \rightarrow 0,
\end{align*}
as $m \rightarrow \infty$. This shows that
$$
\lim_{m \rightarrow \infty}( a \beta_m^l + b \beta_m^s) = \frac{a
- b}{\mu_{k_l}-\mu_{k_s}}M^{ij}.
$$
When $(a,b)=(1,0)$ we obtain Case 5, and when $(a,b)=(0,1)$ we
obtain Case 6. \hfill \qed
\end{description}
\end{proof}

\section{Interactions between tensors and eigenvalues}

The interactions that we will investigate between the types of
tensors defined in Section~\ref{NAPR} and the eigenvalues of
symmetric matrices lead naturally to two families of linear maps.
Each of the next two subsections focuses on one of these families
and explains how it arises.

\subsection{A family of linear maps: divided differences}

 Fix a vector $\mu \in \R^n$. In what follows, the equivalence relation between
numbers from $\N_n$ will be determined by vector $\mu$. We define
$k$ linear maps
$$
T \in T^{k,n} \rightarrow T^{(l)}_{\mbox{\rm \scriptsize out}} \in
T^{k+1,n}, \mbox{ for } l=1,2,...,k.
$$
as follows:
% $(k+1)$-tensors on $\R^n$ denoted
%by $T^{(l)}_{\mbox{\rm \scriptsize out}}$, $l=1,...,k$ as follows.
\begin{equation}
\label{defn-out} \big(T^{(l)}_{\mbox{\rm \scriptsize
out}}\big)^{i_1...i_k i_{k+1}} = \left\{
\begin{array}{ll}
0, & \mbox{ if } i_l \sim i_{k+1}, \\
\displaystyle\frac{T^{i_1...i_{l-1}i_{k+1}i_{l+1}...i_k} -
T^{i_1...i_{l-1}i_li_{l+1}...i_k}}{\mu_{i_{k+1}} - \mu_{i_l}}, &
\mbox{ if } i_l \not\sim i_{k+1}.
\end{array}
\right.
\end{equation}
 Notice that if $T$ is a $\mu$-block-constant tensor, then so is
$T^{(l)}_{\mbox{\rm \scriptsize out}}$ for each $l=1,...,k$. The
easy-to-check claim that these maps are linear means that for any
two tensors $T_1, T_2 \in T^{k,n}$ and $\alpha, \beta \in \R$ we
have
\begin{equation}
\label{sumlab} (\alpha T_1 +  \beta T_2)^{(l)}_{\mbox{\rm
\scriptsize out}} = \alpha (T_1)^{(l)}_{\mbox{\rm \scriptsize
out}} + \beta (T_2)^{(l)}_{\mbox{\rm \scriptsize out}}, \,\,\,
\mbox{ for all } l=1,...,k.
\end{equation}

One can of course iterate this definition: on the space
$T^{k+1,n}$ we can define $k+1$ liner maps into $T^{k+2,n}$, and
so on. We will need a way to keep track of that chain process
somehow. A good enumerating tool for our needs turns out to be the
set of permutations on $\N_k$, $\N_{k+1}$,....

Given a permutation $\sigma$ on $\N_k$ we can naturally view it as
a permutation on $\N_{k+1}$ fixing the last element. Let $\tau_l$
 be the transposition $(l,k+1)$, for all $l=1,...,k,k+1$.
 Define $k+1$ permutations, $\sigma_{(l)}$, on $\N_{k+1}$, as follows:
\begin{equation}
\label{permdefn} \sigma_{(l)} = \sigma \tau_l, \mbox{ for }
l=1,...,k,k+1.
\end{equation}
Informally speaking, given the cycle decomposition of $\sigma$, we
obtain $\sigma_{(l)}$, for each $l=1,...,k$, by inserting the
element $k+1$ immediately after the element $l$, and when $l=k+1$,
the permutation $\sigma_{(k+1)}$ fixes the element $k+1$. Clearly
$\sigma_{(l)}^{-1}(k+1)=l$ for all $l$, and
$$
\{ \mbox{All permutations on } \N_{k+1} \} = \{\sigma \tau_l \, |
\, \sigma \mbox{ is a permutation on } \N_k, \,\,\,
l=1,...,k,k+1\}.
$$

\begin{theorem}
\label{dec14b} Let $\{M_m \}$ be a sequence of symmetric matrices
converging to 0, such that $M_m/\|M_m \|$ converges to $M$. Let
$\mu$ be in $\R^n_{\downarrow}$ and $U_m \rightarrow U \in O^n$ be
a sequence of orthogonal matrices such that
\begin{equation*}
\mbox{\rm Diag}\, \mu + M_m = U_m \big( \mbox{\rm Diag}\,
\lambda(\mbox{\rm Diag}\, \mu + M_m) \big) U_m^T, \, \, \mbox{ for
all } \, \, m=1,2,....
\end{equation*}
Then for every block-constant $k$-tensor $T$ on $\R^n$, any
matrices $H_{1}$,...,$H_{k}$, and any permutation $\sigma$ on
$\N_k$ we have
\begin{equation}
\label{kr} \lim_{m \rightarrow \infty} \Big(\frac{U_m
(\Diag^{\sigma} T) U_m^T - \Diag^{\sigma} T
}{\|M_m\|}\Big)[H_{1},...,H_{k}] = \sum_{l=1}^k
\big(\Diag^{\sigma_{(l)}} T^{(l)}_{\mbox{\rm \scriptsize
out}}\big)[H_{1},..., H_{k}, M_{\mbox{\rm \scriptsize out}}].
\end{equation}
\end{theorem}

\begin{proof}
Both sides of Equation~(\ref{kr}) are linear in each argument
$H_s$. That is why it is enough to prove the result when $H_s$,
for $s=1,...,k$, is an arbitrary matrix, $H_{i_sj_s}$, from the
standard basis on $M^n$. Since we have that
\begin{align*}
\Big(U_m (\Diag^{\sigma} T) U_m^T - \Diag^{\sigma} T
\Big)[H_{i_1j_1},...,H_{i_kj_k}] &= (U_m (\Diag^{\sigma} T) U_m^T)^{\substack{i_1...i_k \\
j_1...j_k}} - (\Diag^{\sigma} T)^{\substack{i_1...i_k \\
j_1...j_k}},
\end{align*}
we begin by developing the first term on the right-hand side.

By the definition of the conjugate action and the fact that $T$ is
block-constant, we have
\begin{align*}
(U_m (\Diag^{\sigma} T) U_m^T)^{\substack{i_1...i_k \\
j_1...j_k}} &= \displaystyle\sum_{\substack{p_{\eta},q_{\eta}
= 1 \\ \eta = 1,...,k } }^{n,...,n} (\Diag^{\sigma} T)^{\substack{p_1...p_k \\
q_1...q_k}} \prod_{\nu =1 }^k
U^{i_{\nu}p_{\nu}}_mU^{j_{\nu}q_{\nu}}_m \\
&=\displaystyle\sum_{\substack{p_{\eta} = 1 \\ \eta = 1,...,k }
}^{n,...,n} T^{p_1...p_k} \prod_{\nu =1 }^k
U^{i_{\nu}p_{\nu}}_mU^{j_{\nu}p_{\sigma^{-1}(\nu)}}_m \\
&=\displaystyle\sum_{\substack{p_{\eta} = 1 \\ \eta = 1,...,k }
}^{n,...,n} T^{p_1...p_k} \prod_{\nu =1 }^k
U^{i_{\nu}p_{\nu}}_mU^{j_{\sigma(\nu)}p_{\nu}}_m \\
&= \displaystyle\sum_{\substack{t_{\eta}=1 \\
\eta=1,...,k}}^{r,...,r} T^{\iota_{t_1}...\iota_{t_k}} \prod_{\nu
=1 }^k \Big( \sum_{p_{\nu} \in I_{t_{\nu}}}
U^{i_{\nu}p_{\nu}}_mU^{j_{\sigma(\nu)}p_{\nu}}_m \Big).
\end{align*}
Thus, we have to take the limit as $m$ approaches infinity of the
expression:
\begin{align*}
\frac{(U_m (\Diag^{\sigma} T) U_m^T - \Diag^{\sigma} T)^{\substack{i_1...i_k \\
j_1...j_k}} }{\|M_m\|} =
\frac{\displaystyle\sum_{t_1,...,t_k=1}^{r,...,r}
T^{\iota_{t_1}...\iota_{t_k}} \prod_{\nu =1 }^k \Big(
\sum_{p_{\nu} \in I_{t_{\nu}}}
U^{i_{\nu}p_{\nu}}_mU^{j_{\sigma(\nu)}p_{\nu}}_m \Big) -
(\Diag^{\sigma} T)^{\substack{i_1...i_k \\ j_1...j_k}} }{\|M_m\|}.
\end{align*}

Assume that $i_l \in I_{v_l}$ and $j_{\sigma(l)} \in I_{s_l}$ for
all $l=1,...,k$.

We investigate several possibilities. Suppose first that among the
pairs
\begin{equation}
\label{pairs:0}
(i_1,j_{\sigma(1)}),(i_2,j_{\sigma(2)}),...,(i_k,j_{\sigma(k)})
\end{equation}
at least two have nonequal entries. It will become clear, that
without loss of generality we may assume they are
$(i_1,j_{\sigma(1)})$ and $(i_2,j_{\sigma(2)})$, that is, $i_1
\not= j_{\sigma(1)}$ and $i_2 \not= j_{\sigma(2)}$. Using
Expansion~(\ref{exp}), for any $t_1$, $t_2$ we observe that:
\begin{align*}
&\lim_{m \rightarrow \infty}\frac{1}{\|M_m\|} \Big(\sum_{p_{1} \in
I_{t_{1}}}^{} U^{i_{1}p_{1}}_mU^{j_{\sigma(1)}p_{1}}_m
\Big)\Big(\sum_{p_{2} \in I_{t_{2}}}^{}
U^{i_{2}p_{2}}_mU^{j_{\sigma(2)}p_{2}}_m \Big) \\
& \hspace{0.5cm} =\lim_{m \rightarrow \infty} \frac{1}{\|M_m\|}
\Big(\delta_{i_1j_{\sigma(1)}} \delta_{v_1t_1} +
\frac{\delta_{v_1t_1} - \delta_{s_1t_1}}{\mu_{i_1} -
\mu_{j_{\sigma(1)}}}
M^{i_1j_{\sigma(1)}} \|M_m\| +o(\|M_m\|) \Big) \times \\
& \hspace{8cm} \Big(\delta_{i_2j_{\sigma(2)}} \delta_{v_2t_2} +
\frac{\delta_{v_2t_2} - \delta_{s_2t_2}}{\mu_{i_2} -
\mu_{j_{\sigma(2)}}} M^{i_2j_{\sigma(2)}} \|M_m\| +o(\|M_m\|)\Big) \\
& \hspace{0.5cm} =\lim_{m \rightarrow \infty} \frac{1}{\|M_m\|}
\Big( \frac{\delta_{v_1t_1} - \delta_{s_1t_1}}{\mu_{i_1} -
\mu_{j_{\sigma(1)}}} M^{i_1j_{\sigma(1)}} \|M_m\| +o(\|M_m\|)
\Big) \Big( \frac{\delta_{v_2t_2} - \delta_{s_2t_2}}{\mu_{i_2} -
\mu_{j_{\sigma(2)}}} M^{i_2j_{\sigma(2)}} \|M_m\| +o(\|M_m\|)\Big) \\
& \hspace{0.5cm} =0.
\end{align*}
Since in this case by definition $(\Diag^{\sigma}
T)^{\substack{i_1...i_k \\ j_1...j_k}}= 0$ we see that the whole
limit above in zero.

Suppose now, that exactly one pair has unequal entries and let it
be $(i_l,j_{\sigma(l)})$. We consider two subcases depending on
whether or not $i_l$ and $j_{\sigma(l)}$ are in the same block.

If both $i_l$ and $j_{\sigma(l)}$ are in one block, that is
$v_l=s_l$, then using Expansion~(\ref{exp}), for arbitrary $t$, we
obtain:
\begin{align*}
\lim_{m \rightarrow \infty}\frac{1}{\|M_m\|} \Big(\sum_{p \in
I_{t}}^{} U^{i_{l}p}_mU^{j_{\sigma(l)}p}_m \Big) &= \lim_{m
\rightarrow \infty}\frac{1}{\|M_m\|} \Big(
\delta_{i_lj_{\sigma(l)}} \delta_{v_lt} + \frac{\delta_{v_lt} -
\delta_{s_lt}}{\mu_{i_l} - \mu_{j_{\sigma(l)}}}
M^{i_lj_{\sigma(l)}} \|M_m\| +o(\|M_m\|)\Big) \\
&= \lim_{m \rightarrow \infty}\frac{o(\|M_m\|)}{\|M_m\|} \\
&=0.
\end{align*}
In this subcase we again have $(\Diag^{\sigma} T)^{\substack{i_1...i_k \\
j_1...j_k}}= 0$, thus the whole limit above is zero.

If $i_l$ and $j_{\sigma(l)}$ are in different blocks, $v_l \not=
s_l$, then $(\Diag^{\sigma} T)^{\substack{i_1...i_k \\
j_1...j_k}}=0$ and by Expansion~(\ref{exp}) we obtain:
\begin{align*}
& \lim_{m \rightarrow \infty} \frac{1}{\|M_m\|} \Big(
\displaystyle\sum_{t_1,...,t_k=1}^{r,...,r}
T^{\iota_{t_1}...\iota_{t_k}} \prod_{\nu =1 }^k \Big(
\sum_{p_{\nu} \in I_{t_{\nu}}}
U^{i_{\nu}p_{\nu}}_mU^{j_{\sigma(\nu)}p_{\nu}}_m
\Big)\Big) \\
& \hspace{1.2cm} =\lim_{m \rightarrow \infty}
\frac{1}{\|M_m\|}\Big( \displaystyle\sum_{t_1,...,t_k=1}^{r,...,r}
T^{\iota_{t_1}...\iota_{t_k}} \prod_{\nu =1 }^k
\Big(\delta_{i_{\nu}j_{\sigma(\nu)}} \delta_{v_{\nu}t_{\nu}} +
\frac{\delta_{v_{\nu}t_{\nu}} -
\delta_{s_{\nu}t_{\nu}}}{\mu_{i_{\nu}} - \mu_{j_{\sigma(\nu)}}}
M^{i_{\nu}j_{\sigma(\nu)}} \|M_m\| +o(\|M_m\|) \Big)\Big).
\end{align*}
We show that the limit of at most two terms in the above sum are
non-zero. Indeed, summands corresponding to $k$-tuples
$(t_1,...,t_k)$ with $t_l \not\in \{v_l,s_l\}$ are equal to zero,
because $\delta_{i_{l}j_{\sigma(l)}}=0$, $\delta_{v_{l}t_{l}}
=\delta_{s_{l}t_{l}} = 0$, and therefore
$$
\delta_{i_{l}j_{\sigma(l)}} \delta_{v_{l}t_{l}} +
\frac{\delta_{v_{l}t_{l}} - \delta_{s_{l}t_{l}}}{\mu_{i_{l}} -
\mu_{j_{\sigma(l)}}} M^{i_{l}j_{\sigma(l)}} \|M_m\| +o(\|M_m\|) =
o(\|M_m\|).
$$

\noindent Similarly, summands corresponding to $k$-tuples
$(t_1,...,t_k)$ with $t_{\nu} \not= v_{\nu}$ for some $\nu \not=
l$ are equal to zero, since then $\delta_{v_{\nu}t_{\nu}} =
\delta_{s_{\nu}t_{\nu}} = 0$ (recall that $v_{\nu} = s_{\nu}$ for
all $\nu \not= l$). Thus, the summands with possible nonzero limit
correspond to the $k$-tuples
$(v_1,...,v_{l-1},v_l,v_{l+1},...,v_k)$ and
$(v_1,...,v_{l-1},s_l,v_{l+1},...,v_k)$. On the other hand, if
$t_{\nu} = v_{\nu} (=s_{\nu})$ for some $\nu \not= l$, then
\begin{equation*}
\delta_{i_{\nu}j_{\sigma(\nu)}} \delta_{v_{\nu}t_{\nu}} +
\frac{\delta_{v_{\nu}t_{\nu}} -
\delta_{s_{\nu}t_{\nu}}}{\mu_{i_{\nu}} - \mu_{j_{\sigma(\nu)}}}
M^{i_{\nu}j_{\sigma(\nu)}} \|M_m\| +o(\|M_m\|) = 1 + o(\|M_m\|).
\end{equation*}
Thus, we can calculate that the above limit is equal to
\begin{align*}
%&\lim_{m \rightarrow \infty} \frac{1}{\|M_m\|} \Big(
%\Big(\frac{T^{k_{v_1}...k_{v_{l-1}}k_{v_l}k_{t_{v+1}}...k_{v_k}}}{\mu_{i_{l}}
%- \mu_{j_{l}}} M^{i_{l}j_{l}} \|M_m\| +o(\|M_m\|)\Big) -
%\Big(\frac{T^{k_{v_1}...k_{v_{l-1}}k_{s_l}k_{t_{v+1}}...k_{v_k}}}{\mu_{i_{l}}
%- \mu_{j_{l}}} M^{i_{l}j_{l}} \|M_m\| +o(\|M_m\|)\Big)\Big) \\
%&=
\frac{T^{\iota_{v_1}...\iota_{v_{l-1}}\iota_{v_l}\iota_{v_{l+1}}...\iota_{v_k}}
-T^{\iota_{v_1}...\iota_{v_{l-1}}\iota_{s_l}\iota_{v_{l+1}}...\iota_{v_k}}
}{\mu_{i_{l}} - \mu_{j_{\sigma(l)}}}M^{i_{l}j_{\sigma(l)}} &=
\frac{T^{i_1...i_{l-1}i_li_{l+1}...i_k} -
T^{i_1...i_{l-1}j_{\sigma(l)}i_{l+1}...i_k} }{\mu_{i_{l}} -
\mu_{j_{\sigma(l)}}}M^{i_{l}j_{\sigma(l)}} \\
&= \frac{T^{i_1...i_{l-1}i_li_{l+1}...i_k} -
T^{i_1...i_{l-1}j_{\sigma(l)}i_{l+1}...i_k} }{\mu_{i_{l}} -
\mu_{j_{\sigma(l)}}}M^{i_{l}j_{\sigma(l)}}_{\mbox{\rm \scriptsize
out}}
\end{align*}
where the first equality follows from the block-constant structure
of $T$ and the second from the premise in this case that $i_l$ and
$j_{\sigma(l)}$ are in different blocks.

In the last case when $i_{\nu} = j_{\sigma(\nu)}$ for all $\nu =
1,...,k$,
 the limit is equal to
\begin{align*}
\lim_{m \rightarrow \infty} \frac{1}{\|M_m\|} \Big(
 T^{i_1...i_k} (1+o(\|M_m\|)) - T^{i_1...i_k}
\Big) = 0.
\end{align*}
\noindent With that we finished calculating the limit in the
left-hand side of Equation~(\ref{kr}).

We now compute the right-hand side of Equation~(\ref{kr}) and
compare with the results above. Suppose that $\sigma(l)=m$, then
by the definition of $\sigma_{(l)}$ we have
$\sigma^{-1}_{(l)}(m)=k+1$, $\sigma^{-1}_{(l)}(k+1)=l$, and for
any integer $i \in \N_{k+1} \backslash \{k+1,m \}$ we have
$\sigma^{-1}_{(l)}(i)= \sigma^{-1}(i)$. Below we use the standard
notation that a ``hat'' above a term in a product means that the
term is omitted. Since $\sigma_{(l)}^{-1}(k+1)=l \not= k+1$ we use
the second part of Lemma~\ref{dot-prod} to compute:
\begin{align*}
\sum_{l=1}^k  \big(\Diag^{\sigma_{(l)}} T^{(l)}_{\mbox{\rm
\scriptsize out}}\big) [H_{i_1j_1}&,...,H_{i_kj_k}, M_{\mbox{\rm
\scriptsize out}}] =\sum_{l=1}^k \langle T^{(l)}_{\mbox{\rm
\scriptsize out}}, H_{i_1j_1} \circ_{\sigma_{(l)}} \cdots
\circ_{\sigma_{(l)}} H_{i_kj_k} \circ_{\sigma_{(l)}} M_{\mbox{\rm
\scriptsize out}}
\rangle \\
&=\sum_{l=1}^k \big(T^{(l)}_{\mbox{\rm \scriptsize
out}}\big)^{i_1...i_kj_{\sigma_{(l)}(k+1)}}
\big(\delta_{i_{1}}j_{\sigma_{(l)}(1)} \cdots
\widehat{\delta_{i_{l}}j_{\sigma_{(l)}(l)}} \cdots
\delta_{i_{k}}j_{\sigma_{(l)}(k)} \big)
M^{j_{\sigma_{(l)}(k+1)}i_{\sigma^{-1}_{(l)}(k+1)}}_{\mbox{\rm
\scriptsize out}} \\
&=\sum_{l=1}^k \big(T^{(l)}_{\mbox{\rm \scriptsize
out}}\big)^{i_1...i_kj_{\sigma(l)}}
\big(\delta_{i_{1}j_{\sigma(1)}} \cdots
\widehat{\delta_{i_{l}j_{\sigma(l)}}} \cdots
\delta_{i_{k}j_{\sigma(k)}} \big) M^{j_{\sigma(l)}i_l}_{\mbox{\rm
\scriptsize out}}.
\end{align*}
It is clear that if at least two of the pairs
$(i_1,j_{\sigma(1)}),(i_2,j_{\sigma(2)}),...,(i_k,j_{\sigma(k)})$
have different entries, then the sum is zero. Let now exactly on
of the pairs have unequal entries, say $i_l \not= j_{\sigma(l)}$,
then the above sum will be equal to
$$
\big(T^{(l)}_{\mbox{\rm \scriptsize
out}}\big)^{i_1...i_kj_{\sigma(l)}}
\big(\delta_{i_{1}j_{\sigma(1)}} \cdots
\widehat{\delta_{i_{l}j_{\sigma(l)}}} \cdots
\delta_{i_{k}j_{\sigma(k)}} \big) M^{j_{\sigma(l)}i_l}_{\mbox{\rm
\scriptsize out}}.
$$
If $i_l$ and $j_{\sigma(l)}$ are in the same block, then
$\big(T^{(l)}_{\mbox{\rm \scriptsize
out}}\big)^{i_1...i_kj_{\sigma(l)}}=0$ by the definition of
$T^{(l)}_{\mbox{\rm \scriptsize out}}$. If $i_l$ and
$j_{\sigma(l)}$ are not in the same block, then the last
expression above is equal to
$$
\big(T^{(l)}_{\mbox{\rm \scriptsize
out}}\big)^{i_1...i_kj_{\sigma(l)}}
M^{j_{\sigma(l)}i_l}_{\mbox{\rm \scriptsize out}} =
\frac{T^{i_1...i_{l-1}i_li_{l+1}...i_k} -
T^{i_1...i_{l-1}j_{\sigma(l)}i_{l+1}...i_k} }{\mu_{i_{l}} -
\mu_{j_{\sigma(l)}}}M^{i_{l}j_{\sigma(l)}}_{\mbox{\rm \scriptsize
out}},
$$
because $M$ is a symmetric matrix. Finally, if
$i_{\nu}=j_{\sigma(\nu)}$ for all $\nu = 1,...,k$, then again
$\big(T^{(l)}_{\mbox{\rm \scriptsize
out}}\big)^{i_1...i_kj_{\sigma(l)}}=0$. These outcomes are equal
to the results in the corresponding cases in the first part of the
proof, the theorem follows.  \hfill \qed
\end{proof}

\subsection{A second family of linear maps: ``inflating'' diagonal hyper-planes}

Recall that $\tau_l$ denotes the transposition $(l,k+1)$ on
$\N_{k+1}$. Fix a vector $\mu \in \R^n$ defining the equivalence
relation on $\N_n$. We define $k$ linear maps
$$
T \in T^{k,n} \rightarrow T^{(l)}_{\mbox{\rm \scriptsize in}} \in
T^{k+1,n}, \mbox{ for } l=1,2,...,k.
$$
%Given a $k$-tensor $T \in T^{k,n}$ and a number $l \in \N_k$, the
%tensor $T^{(l)}_{\mbox{\rm \scriptsize in}}$ is given by the
%formula:
as follows:
\begin{equation}
\label{dec14c:eqn} \big(T^{\tau_l}_{\mbox{\rm \scriptsize
in}}\big)^{i_1...i_ki_{k+1}} = \left\{
\begin{array}{ll}
T^{i_1...i_{l-1}i_{k+1}i_{l+1}...i_k}, & \mbox{ if } i_l \sim
i_{k+1}, \\
0, & \mbox{ if } i_l \not\sim i_{k+1}.
\end{array}
\right.
\end{equation}
Notice that if $T$ is a block-constant tensor, then so is
$T^{\tau_l}_{\mbox{\rm \scriptsize in}}$ for each $l=1,...,k$. In
that case, we clearly have
$$
\big(T^{\tau_l}_{\mbox{\rm \scriptsize
in}}\big)^{i_1...i_ki_{k+1}} = \left\{
\begin{array}{ll}
T^{i_1...i_{l-1}i_{l}i_{l+1}...i_k}, & \mbox{ if } i_l \sim
i_{k+1}, \\
0, & \mbox{ if } i_l \not\sim i_{k+1}.
\end{array}
\right.
$$
It is again easy to check that these maps are linear, that is, for
any two tensors $T_1, T_2 \in T^{k,n}$ and $\alpha, \beta \in \R$
we have
\begin{equation*}
\label{sumlab-a} (\alpha T_1 +  \beta T_2)^{(l)}_{\mbox{\rm
\scriptsize in}} = \alpha (T_1)^{(l)}_{\mbox{\rm \scriptsize in}}
+ \beta (T_2)^{(l)}_{\mbox{\rm \scriptsize in}}, \,\,\, \mbox{ for
all } l=1,...,k.
\end{equation*}
Define also
\begin{equation}
\label{perm-lift} \big(T^{\tau_l}\big)^{i_1...i_ki_{k+1}} =
\left\{
\begin{array}{ll}
T^{i_1...i_{l-1}i_{l}i_{l+1}...i_k}, & \mbox{ if } i_l =
i_{k+1}, \\
0, & \mbox{ if } i_l \not= i_{k+1}.
\end{array}
\right.
\end{equation}
In other words, $T^{\tau_l}$ is a $(k+1)$-tensor with entries off
the hyper plane $i_l=i_{k+1}$ equal to zero. On the hyper plane
$i_l=i_{k+1}$ we have placed the original tensor $T$.

Before we formulate the main result or this subsection we need two
technical lemmas.

\begin{proposition}
\label{prop11} Let $T$ be any $k+1$-tensor, $x$ be any vector in
$\R^n$, let $V$ be any orthogonal matrix, and $\sigma$ a
permutation on $\N_k$. Then the following identity holds:
$$
V\big(\Diag^{\sigma} (T[x])\big)V^T = \big(
V(\Diag^{\sigma_{(k+1)}} T)V^T \big)[V(\Diag x) V^T].
$$
\end{proposition}

\begin{proof}
Let $H_{i_1j_1}$,...,$H_{i_kj_k}$ be any $k$ basic matrices.
Recall that $\sigma_{(k+1)}(i)=\sigma(i)$ for all $i \in \N_{k}$
and $\sigma_{(k+1)}(k+1)=k+1$.  Using Theorem~\ref{jen3a} twice,
we compute
\begin{align*}
\big(V\big(\Diag^{\sigma} &(T[x])\big)V^T\big)^{\substack{i_1...i_k \\
j_1...j_k}} = \big(V\big(\Diag^{\sigma}
(T[x])\big)V^T\big)[H_{i_1j_1},...,H_{i_kj_k}] \\
&=\langle T[x], \tilde{H}_{i_1j_1} \circ_{\sigma}...\circ_{\sigma} \tilde{H}_{i_kj_k} \rangle \\
&= \sum_{p_1,...,p_k=1}^{n,...,n} (T[x])^{p_1...p_k}
\tilde{H}_{i_1j_1}^{p_1p_{\sigma^{-1}(1)}} \cdots
\tilde{H}_{i_kj_k}^{p_kp_{\sigma^{-1}(k)}} \\
&=\sum_{p_1,...,p_k,p_{k+1}=1}^{n,...,n} T^{p_1...p_{k+1}}
x^{p_{k+1}} \tilde{H}_{i_1j_1}^{p_1p_{\sigma^{-1}(1)}} \cdots
\tilde{H}_{i_kj_k}^{p_kp_{\sigma^{-1}(k)}} \\
&=\sum_{p_1,...,p_k,p_{k+1}=1}^{n,...,n} T^{p_1...p_{k+1}}
 \tilde{H}_{i_1j_1}^{p_1p_{\sigma_{(k+1)}^{-1}(1)}} \cdots
\tilde{H}_{i_kj_k}^{p_kp_{\sigma_{(k+1)}^{-1}(k)}}
(\Diag x)^{p_{k+1}p_{\sigma_{(k+1)}^{-1}(k+1)}}\\
&=\langle T, \tilde{H}_{i_1j_1} \circ_{\sigma_{(k+1)}}...\circ_{\sigma_{(k+1)}}
\tilde{H}_{i_kj_k} \circ_{\sigma_{(k+1)}} \Diag x \rangle\\
&=\big(V (\Diag^{\sigma_{(k+1)}} T)
V^T\big)[H_{i_1j_1},...,H_{i_kj_k}, V(\Diag x)V^T] \\
&=\big(\big( V(\Diag^{\sigma_{(k+1)}} T)V^T \big) [V(\Diag x)
V^T]\big)^{\substack{i_1...i_k \\ j_1...j_k}}.
\end{align*}
Since these equalities hold for all $i_1...i_k$ and $j_1...j_k$ we
are done. \hfill \qed
\end{proof}

The next lemma says that for any block-constant tensor $T$,
$\Diag^{\sigma} T$ is invariant under conjugations with a
block-diagonal orthogonal matrix.

\begin{lemma}
\label{invarLem} Let $T$ be a block-constant $k$-tensor on $\R^n$,
let $U \in O^n$ be a block-diagonal matrix (both with respect to
one and the same partitioning of $\N_n$). Then for any permutation
$\sigma$ in $\N_k$ we have the identity
$$
U(\Diag^{\sigma} T)U^T = \Diag^{\sigma} T.
$$
\end{lemma}

\begin{proof}
Let $\{I_1$,...,$I_r\}$ be the partitioning of the integers $\N_n$
that determines the block structure. Notice that $U^{ip}U^{jp}=0$
whenever $i \not\sim j$ or $i \not\sim p$, and that $\sum_{p \in
I_s}U^{ip}U^{jp} = \delta_{ij}$ whenever $i \in I_s$. Let
$(i_1,...,i_k)$ be an arbitrary multi index and suppose that $i_s
\in I_{\nu_s}$ for $s=1,...,k$.  We expand the left-hand side of
the identity:
\begin{align*}
\big(U(\Diag^{\sigma} T)U^T\big)^{\substack{i_1...i_k \\
j_1...j_k}}  &= \sum_{\substack{p_s, q_s = 1 \\
s=1,...,k}}^{n,...,n} (\Diag^{\sigma} T)^{\substack{p_1...p_k \\
q_1...q_k}} U^{i_1p_1}U^{j_1q_1} \cdots U^{i_kp_k}U^{j_kq_k} \\
&= \sum_{p_1,...,p_k = 1}^{n,...,n} T^{p_1...p_k} U^{i_1p_1}U^{j_1p_{\sigma^{-1}(1)}} \cdots U^{i_kp_k}U^{j_kp_{\sigma^{-1}(k)}} \\
&= \sum_{p_1,...,p_k = 1}^{n,...,n} T^{p_1...p_k} U^{i_1p_1}U^{j_{\sigma(1)}p_{1}} \cdots U^{i_kp_k}U^{j_{\sigma(k)}p_{k}} \\
&= \sum_{t_1,...,t_k =1}^{r,...,r} T^{\iota_{t_1}...\iota_{t_k}}
\sum_{\substack{p_l \in I_{t_l}
\\ l=1,...,k}} U^{i_1p_1}U^{j_{\sigma(1)}p_{1}} \cdots U^{i_kp_k}U^{j_{\sigma(k)}p_{k}} \\
&= T^{\iota_{\nu_1}...\iota_{\nu_k}} \sum_{\substack{p_l \in
I_{\nu_l} \\ l=1,...,k}} U^{i_1p_1}U^{j_{\sigma(1)}p_{1}} \cdots U^{i_kp_k}U^{j_{\sigma(k)}p_{k}} \\
&= T^{\iota_{\nu_1}...\iota_{\nu_k}} \delta_{i_1j_{\sigma(1)}}
\cdots \delta_{i_kj_{\sigma(k)}} \\
&= T^{i_1...i_k} \delta_{i_1j_{\sigma(1)}} \cdots \delta_{i_kj_{\sigma(k)}} \\
&= \big(\Diag^{\sigma} T \big)^{\substack{i_1...i_k \\
j_1...j_k}}.
\end{align*}
The next to the last equality follows from the fact that $T$ is
block-constant. Since the multi index $(i_1,...,i_k,j_1,...,j_k)$
was arbitrary, the claim in the lemma follows. \hfill \qed
\end{proof}

\begin{theorem}
\label{dec15a}
%Let $\{M_m \}$ be a sequence of symmetric matrices converging to
%0, such that $M_m/\|M_m \|$ converges to $M$. Let $\mu$ be in
%$\R^n_{\downarrow}$ and $U_m \rightarrow U \in O^n$ be a sequence
%of orthogonal matrices such that
%\begin{equation*}
%\mbox{\rm Diag}\, \mu + M_m = U_m \big( \mbox{\rm Diag}\,
%\lambda(\mbox{\rm Diag}\, \mu + M_m) \big) U_m^T, \, \, \mbox{ for
%all } \, \, m=1,2,....
%\end{equation*}
%Let $h$ be the vector defined by (\ref{limh}),
Let $U \in O^n$ be a block-diagonal orthogonal matrix. Let $M$ be
an arbitrary symmetric matrix, and let $h \in \R^n$ be a vector,
such that
\begin{equation}
\label{bg} U^TM_{\mbox{\rm \scriptsize in}}U = \Diag h.
\end{equation}
Let $H_{1}$,...,$H_{k}$ be arbitrary matrices, and let $\sigma$ be
a permutation on $\N_k$. Then
\begin{enumerate}
\item for any block-constant $(k+1)$-tensor $T$ on $\R^n$,
$$
\langle T[h], \tilde{H}_1 \circ_{\sigma} \cdots \circ_{\sigma}
\tilde{H}_k \rangle = \langle T, H_1 \circ_{\sigma_{(k+1)}} \cdots
\circ_{\sigma_{(k+1)}} H_k \circ_{\sigma_{(k+1)}} M_{\mbox{\rm
\scriptsize in}} \rangle
$$
\item for any block-constant $k$-tensor $T$ on $\R^n$
$$
\langle T^{\tau_l}[h], \tilde{H}_1 \circ_{\sigma} \cdots
\circ_{\sigma} \tilde{H}_k \rangle  = \langle
T^{\tau_l}_{\mbox{\rm \scriptsize in}}, H_1 \circ_{\sigma_{(l)}}
\cdots \circ_{\sigma_{(l)}} H_k \circ_{\sigma_{(l)}} M_{\mbox{\rm
\scriptsize in}} \rangle, \,\,\, \mbox{ for all } l=1,...,k,
$$
\end{enumerate}
where the permutations $\sigma_{(l)}$, for $l=0,1,...,k$ are
defined by (\ref{permdefn}), $\tilde{H}_i=U^TH_iU$ for all
$i=1,...,k$.
%and $\tilde{M}_{\mbox{\rm \scriptsize in}} =
%U^TM_{\mbox{\rm \scriptsize in}}U$.
\end{theorem}

\begin{proof}
To see that the first identity holds we use Theorem~\ref{jen3a},
Proposition~\ref{prop11}, Formula~(\ref{bg}), and
Lemma~\ref{invarLem} in that order, as follows:
\begin{align*}
\langle T[h], \tilde{H}_1 \circ_{\sigma} \cdots \circ_{\sigma}
\tilde{H}_k \rangle &= \big(U(\Diag^{\sigma} T[h])U^T
\big)[H_1,...,H_k] \\
&= \big(U(\Diag^{\sigma_{(k+1)}}T)U^T \big)[H_1,...,H_k, U(\Diag
h)U^T] \\
&= \big(U(\Diag^{\sigma_{(k+1)}}T)U^T \big)[H_1,...,H_k,
M_{\mbox{\rm \scriptsize in}}] \\
&= \big(\Diag^{\sigma_{(k+1)}}T \big)[H_1,...,H_k,
M_{\mbox{\rm \scriptsize in}}] \\
&= \langle T, H_1 \circ_{\sigma_{(k+1)}} \cdots
\circ_{\sigma_{(k+1)}} H_k \circ_{\sigma_{(k+1)}} M_{\mbox{\rm
\scriptsize in}} \rangle.
\end{align*}
The last equality follows again from Theorem~\ref{jen3a}.

To show the second identity, since both sides are linear in $H_s$
for every $s=1,...,k$, it is enough to prove it only in the case
when $H_s$ is an arbitrary basic matrix $H_{i_sj_s}$. Fix $k$
basic matrices $H_{i_1j_1}$,...,$H_{i_kj_k}$ and suppose that $i_s
\in I_{\nu_s}$ for $s=1,...,k$. The left-hand side is equal to
\begin{align*}
\langle T^{\tau_l}[h], \tilde{H}_{i_1j_1} \circ_{\sigma} \cdots
\circ_{\sigma} \tilde{H}_{i_kj_k} \rangle &= \big(U(\Diag^{\sigma}
T^{\tau_l}[h])U^T\big)[H_{i_1j_1},...,H_{i_kj_k}] \\
&= \big(U(\Diag^{\sigma} T^{\tau_l}[h]
)U^T\big)^{\substack{i_1...i_k \\ j_1...j_k}} \\
&= \sum_{\substack{p_1,...,p_k = 1 \\ q_1,...,q_k = 1}}^{n,...,n}
(\Diag^{\sigma} T^{\tau_l}[h] )^{\substack{p_1...p_k \\
q_1...q_k}}U^{i_1p_1}U^{j_1q_1}
\cdots U^{i_kp_k}U^{j_kq_k} \\
&= \sum_{p_1,...,p_k = 1}^{n,...,n}
(T^{\tau_l}[h])^{p_1...p_k}U^{i_1p_1}U^{j_1p_{\sigma^{-1}(1)}}
\cdots U^{i_kp_k}U^{j_kp_{\sigma^{-1}(k)}} \\
&= \sum_{p_1,...,p_k = 1}^{n,...,n}
(T^{\tau_l}[h])^{p_1...p_k}U^{i_1p_1}U^{j_{\sigma(1)}p_{1}}
\cdots U^{i_kp_k}U^{j_{\sigma(k)}p_{k}} \\
&= \sum_{p_1,...,p_k = 1}^{n,...,n} \sum_{p_{k+1}=1}^n
(T^{\tau_l})^{p_1...p_kp_{k+1}}h^{p_{k+1}}U^{i_1p_1}U^{j_{\sigma(1)}p_{1}}
\cdots U^{i_kp_k}U^{j_{\sigma(k)}p_{k}} \\
&= \sum_{p_1,...,p_k = 1}^{n,...,n}
T^{p_1...p_k}h^{p_{l}}U^{i_1p_1}U^{j_{\sigma(1)}p_{1}}
\cdots U^{i_kp_k}U^{j_{\sigma(k)}p_{k}} \\
&= \sum_{s_1,...,s_k=1}^{r,...,r} T^{\iota_{s_1}...\iota_{s_k}}
\sum_{\substack{p_{\eta} \in I_{s_{\eta}} \\ \eta=1,...,k}}
h^{p_{l}}U^{i_1p_1}U^{j_{\sigma(1)}p_{1}}
\cdots U^{i_kp_k}U^{j_{\sigma(k)}p_{k}} \\
&=  T^{i_1...i_k} \delta_{i_1j_{\sigma(1)}} \cdots
\widehat{\delta_{i_lj_{\sigma(l)}}} \cdots
\delta_{i_kj_{\sigma(k)}} \sum_{p_{l} \in I_{\nu_{l}}}
h^{p_{l}}U^{i_lp_l}U^{j_{\sigma(l)}p_{l}} \\
&= T^{i_1...i_k} \delta_{i_1j_{\sigma(1)}} \cdots
\widehat{\delta_{i_lj_{\sigma(l)}}} \cdots
\delta_{i_kj_{\sigma(k)}} M^{i_lj_{\sigma(l)}}_{\mbox{\rm
\scriptsize in}}.
\end{align*}
Now we evaluate the right-hand side of the identity. We will use
the second part of Lemma~\ref{dot-prod} since
$\sigma_{(l)}^{-1}(k+1)=l \not= k+1$. Recall also that
$\sigma_{(l)}(s)=\sigma(s)$ for $s \in \N_{k+1}\backslash
\{l,k+1\}$ and $\sigma_{(l)}(k+1)=\sigma(l)$ for all $l=1,...,k$.
\begin{align*}
\langle T^{\tau_l}_{\mbox{\rm \scriptsize in}}, H_{i_1j_1}
\circ_{\sigma_{(l)}} \cdots \circ_{\sigma_{(l)}} H_{i_kj_k}
\circ_{\sigma_{(l)}} & M_{\mbox{\rm \scriptsize in}} \rangle \\
&= \big(T^{\tau_l}_{\mbox{\rm \scriptsize
in}}\big)^{i_1...i_kj_{\sigma_{(l)}(k+1)}}
\delta_{i_1j_{\sigma_{(l)}(1)}} \cdots
\widehat{\delta_{i_lj_{\sigma_{(l)}(l)}}}\cdots
\delta_{i_kj_{\sigma_{(l)}(k)}} M_{\mbox{\rm \scriptsize in}}^{j_{\sigma_{(l)}(k+1)}i_{\sigma_{(l)}^{-1}(k+1)}} \\
&= \big(T^{\tau_l}_{\mbox{\rm \scriptsize
in}}\big)^{i_1...i_kj_{\sigma(l)}} \delta_{i_1j_{\sigma(1)}}
\cdots \widehat{\delta_{i_lj_{k+1}}}\cdots
\delta_{i_kj_{\sigma(k)}} M_{\mbox{\rm \scriptsize in}}^{j_{\sigma(l)}i_{l}} \\
&= T^{i_1...i_k} \delta_{i_1j_{\sigma(1)}} \cdots
\widehat{\delta_{i_lj_{k+1}}}\cdots
\delta_{i_kj_{\sigma(k)}} M_{\mbox{\rm \scriptsize in}}^{j_{\sigma(l)}i_{l}} \\
&= T^{i_1...i_k} \delta_{i_1j_{\sigma(1)}} \cdots
\widehat{\delta_{i_lj_{k+1}}}\cdots
\delta_{i_kj_{\sigma(k)}} M_{\mbox{\rm \scriptsize in}}^{i_{l}j_{\sigma(l)}} \\
&= T^{i_1...i_k} \delta_{i_1j_{\sigma(1)}} \cdots
\widehat{\delta_{i_lj_{\sigma(l)}}}\cdots
\delta_{i_kj_{\sigma(k)}} M_{\mbox{\rm \scriptsize
in}}^{i_{l}j_{\sigma(l)}}.
\end{align*}
In the third equality above we used the fact that $T$ is
block-constant, plus the fact that $M_{\mbox{\rm \scriptsize
in}}^{j_{\sigma(l)}i_{l}}=0$ if $j_{\sigma(l)} \not\sim i_{l}$. In
the fourth we used the fact that $M$ is a symmetric matrix. The
last equality holds because we changed the format of the missing
multiple, while keeping the present multiples the same. \hfill
\qed
\end{proof}

\begin{proposition}
\label{dec15b} \rm Let $U \in O(n)$ be an block-diagonal
orthogonal matrix, let $H$ be an arbitrary $n \times n$ matrix,
and $\sigma$ an arbitrary permutation on $\N_k$.
\begin{enumerate}
\item If $T$ is a $(k+1)$-tensor such that for some fixed $l \in
\N_k$ we have $T^{p_1...p_l...p_{k+1}}=0$ whenever $p_l \sim
p_{k+1}$, then
$$
\big(U (\Diag^{\sigma_{(l)}} T) U^T \big)[H_{\mbox{\rm \scriptsize
in}}] = 0.
$$
\item If $T$ is a $(k+1)$-tensor such that for some fixed $l \in
\N_k$ we have $T^{p_1...p_l...p_{k+1}}=0$ whenever $p_l \not\sim
p_{k+1}$, then
$$
\big(U (\Diag^{\sigma_{(l)}} T) U^T \big)[H_{\mbox{\rm \scriptsize
out}}] = 0.
$$
\item If $T$ is any $(k+1)$-tensor, then
$$
\big(U (\Diag^{\sigma_{(k+1)}} T) U^T \big)[H_{\mbox{\rm
\scriptsize out}}] = 0.
$$
\end{enumerate}
\end{proposition}

\begin{proof}
Fix an index $l$ in $\N_k$. Let $H_{i_1j_1}$,...,$H_{i_kj_k}$ be
arbitrary basic matrices, and let $H$ be an arbitrary matrix.
Using the definitions we compute.
\begin{align*}
&\big(U (\Diag^{\sigma_{(l)}} T) U^T
\big)[H_{i_1j_1},...,H_{i_kj_k}, H] = \sum_{i_{k+1},
j_{k+1}=1}^{n,n} \big(U (\Diag^{\sigma_{(l)}} T) U^T
\big)^{\substack{i_1...i_{k+1} \\ j_1...j_{k+1}}} H^{i_{k+1}j_{k+1}} \\
&=\sum_{i_{k+1}, j_{k+1}=1}^{n,n} \sum_{\substack{p_s, q_s=1 \\
s=1,...,k+1}}^{n,...,n} (\Diag^{\sigma_{(l)}}
T)^{\substack{p_1...p_{k+1} \\ q_1...q_{k+1}}}
U^{i_1p_1}U^{j_1q_1} \cdots U^{i_{k+1}p_{k+1}}U^{j_{k+1}q_{k+1}} H^{i_{k+1}j_{k+1}} \\
&= \sum_{i_{k+1}, j_{k+1}=1}^{n,n} \sum_{\substack{p_s =1 \\
s=1,...,k+1}}^{n,...,n} T^{p_1...p_{k+1}}
U^{i_1p_1}U^{j_1p_{\sigma^{-1}_{(l)}(1)}} \cdots
U^{i_{k+1}p_{k+1}}U^{j_{k+1}p_{\sigma^{-1}_{(l)}(k+1)}} H^{i_{k+1}j_{k+1}} \\
&= \sum_{i_{k+1}, j_{k+1}=1}^{n,n} \sum_{\substack{p_s =1 \\
s=1,...,k+1}}^{n,...,n} T^{p_1...p_{k+1}}
U^{i_1p_1}U^{j_{\sigma_{(l)}(1)}p_{1}} \cdots
U^{i_lp_l}U^{j_{\sigma_{(l)}(l)}p_{l}} \cdots
U^{i_{k+1}p_{k+1}}U^{j_{\sigma_{(l)}(k+1)}p_{k+1}} H^{i_{k+1}j_{k+1}} \\
&= \sum_{i_{k+1}, j_{k+1}=1}^{n,n} \sum_{\substack{p_s =1 \\
s=1,...,k+1}}^{n,...,n} T^{p_1...p_{k+1}}
U^{i_1p_1}U^{j_{\sigma(1)}p_{1}} \cdots U^{i_lp_l}U^{j_{k+1}p_{l}}
\cdots U^{i_{k+1}p_{k+1}}U^{j_{\sigma(l)}p_{k+1}}
H^{i_{k+1}j_{k+1}}.
\end{align*}
Suppose now that $T$ is a $(k+1)$-tensor with
$T^{p_1...p_l...p_{k+1}}=0$ whenever $p_l \sim p_{k+1}$ and that
$H=H_{\mbox{\rm \scriptsize in}}$. Then $H^{i_{k+1}j_{k+1}}
\not=0$ implies that $i_{k+1} \sim j_{k+1}$. In that case, by the
fact that $U$ is block diagonal,
$U^{j_{k+1}p_{l}}U^{i_{k+1}p_{k+1}} \not= 0$ implies that $p_l
\sim p_{k+1}$, which implies that $T^{p_1...p_l...p_{k+1}}=0$.
Thus every summand in the double sum above is zero.

In the second case, suppose $T$ is a $(k+1)$-tensor with
$T^{p_1...p_l...p_{k+1}}=0$ whenever $p_l \not\sim p_{k+1}$ and
$H=H_{\mbox{\rm \scriptsize out}}$. Then $H^{i_{k+1}j_{k+1}}
\not=0$ implies that $i_{k+1} \not\sim j_{k+1}$. In that case, by
the fact that $U$ is block diagonal,
$U^{j_{k+1}p_{l}}U^{i_{k+1}p_{k+1}} \not= 0$ implies that $p_l
\not\sim p_{k+1}$, which implies that $T^{p_1...p_l...p_{k+1}}=0$.
The sum is zero.

In the third case, suppose that $T$ is any $(k+1)$-tensor and
$H=H_{\mbox{\rm \scriptsize out}}$. A calculation almost identical
to the above one (it differs only in the last step) shows that
\begin{align*}
\big(U (\Diag^{\sigma_{(k+1)}} T&) U^T
\big)[H_{i_1j_1},...,H_{i_kj_k}, H] = \\
&\sum_{i_{k+1}, j_{k+1}=1}^{n,n} \sum_{\substack{p_s =1 \\
s=1,...,k+1}}^{n,...,n} T^{p_1...p_{k+1}}
U^{i_1p_1}U^{j_{\sigma(1)}p_{1}} \cdots
U^{i_kp_k}U^{j_{\sigma(k)}p_{k}}
U^{i_{k+1}p_{k+1}}U^{j_{k+1}p_{k+1}} H^{i_{k+1}j_{k+1}}.
\end{align*}
Then $H^{i_{k+1}j_{k+1}} \not=0$ implies that $i_{k+1} \not\sim
j_{k+1}$. In that case, by the fact that $U$ is block diagonal,
$U^{j_{k+1}p_{k+1}}U^{i_{k+1}p_{k+1}} = 0$. Again the sum is zero.
 \hfill \qed
\end{proof}

The next result is a consequence of Theorem~\ref{jen3a},
Theorem~\ref{dec15a}, and Proposition~\ref{dec15b} applied with
$U=I$ and using the fact that $M=M_{\mbox{\rm \scriptsize
in}}+M_{\mbox{\rm \scriptsize out}}$.

\begin{corollary}
\label{jan11a} Let $U \in O^n$ be a block-diagonal orthogonal
matrix and let $\sigma$ be a permutation on $\N_k$. Let $M$ be an
arbitrary symmetric matrix, and $h \in \R^n$ be a vector, such
that $U^TM_{\mbox{\rm \scriptsize in}}U = \Diag h$. Then
\begin{enumerate}
\item for any block-constant $(k+1)$-tensor $T$ on $\R^n$,
$$
U\big(\Diag^{\sigma} (T[h]) \big)U^T =
\big(\Diag^{\sigma_{(k+1)}}T \big) [M];
$$
\item for any block-constant $k$-tensor $T$ on $\R^n$
$$
U \big(\Diag^{\sigma} (T^{\tau_l}[h]) \big)U^T =
\big(\Diag^{\sigma_{(l)}} T^{\tau_l}_{\mbox{\rm \scriptsize in}}
\big) [M], \,\,\, \mbox{ for all } l=1,...,k,
$$
\end{enumerate}
where the permutations $\sigma_{(l)}$, for $l \in \N_k$, are
defined by (\ref{permdefn}).
\end{corollary}

Notice that if the vector $\mu$, defining the equivalence relation
on $\N_n$, has distinct coordinates, then every tensor from
$T^{k,n}$ is block-constant and the block-diagonal orthogonal
matrices are precisely the signed identity matrices (those with
plus or minus one on the main diagonal and zeros everywhere else).
In this case we also have $i \sim j$ if, and only if, $i = j$ and
thus $T^{\tau_l}_{\mbox{\rm \scriptsize in}} = T^{\tau_l}$.
Moreover, since Proposition~\ref{dec15b} holds for arbitrary
matrices (symmetric or not), we obtain the next result, valid for
an arbitrary matrix $H$.

\begin{corollary}
\label{jan11abc} Let $\sigma$ be a permutation on $\N_k$ and let
$H$ be an arbitrary matrix. Then
\begin{enumerate}
\item for any $(k+1)$-tensor $T$ on $\R^n$,
$$
\Diag^{\sigma} (T[\diag H]) = \big(\Diag^{\sigma_{(k+1)}}T \big)
[H];
$$
\item for any $k$-tensor $T$ on $\R^n$
$$
\Diag^{\sigma} (T^{\tau_l}[\diag H]) = \big(\Diag^{\sigma_{(l)}}
T^{\tau_l} \big) [H], \,\,\, \mbox{ for all } l=1,...,k,
$$
\end{enumerate}
where the permutations $\sigma_{(l)}$, for $l \in \N_k$, are
defined by (\ref{permdefn}).
\end{corollary}

Finally, combining Theorem~\ref{dec14b} and
Proposition~\ref{dec15b} we get the next corollary.

\begin{corollary}
\label{dec14bc} Let $\{M_m \}$ be a sequence of symmetric matrices
converging to 0, such that $M_m/\|M_m \|$ converges to $M$. Let
$\mu$ be in $\R^n_{\downarrow}$ and $U_m \rightarrow U \in O^n$ be
a sequence of orthogonal matrices such that
\begin{equation*}
\mbox{\rm Diag}\, \mu + M_m = U_m \big( \mbox{\rm Diag}\,
\lambda(\mbox{\rm Diag}\, \mu + M_m) \big) U_m^T, \, \, \mbox{ for
all } \, \, m=1,2,....
\end{equation*}
Then for every block-constant $k$-tensor $T$ on $\R^n$, and any
permutation $\sigma$ on $\N_k$ we have
\begin{equation}
\label{kr} \lim_{m \rightarrow \infty} \frac{U_m (\Diag^{\sigma}
T) U_m^T - \Diag^{\sigma} T }{\|M_m\|} = \sum_{l=1}^k
\big(\Diag^{\sigma_{(l)}} T^{(l)}_{\mbox{\rm \scriptsize
out}}\big)[M].
\end{equation}
\end{corollary}

\section{A determinant connection}

Let $T \in T^{k,n}$ be any $\mu$-symmetric tensor on $R^n$. A
little bit of thought shows that $T$ can be decomposed into $T = A
+ B$, where $A,B \in T^{k,n}$, $A$ is a block-constant tensor and
$B$ has the following property: for every multi index
$(i_1,...,i_k)$ with pairwise distinct entries, $B^{i_1...i_k} =
0$. In this section we investigate a rather curious fact about
tensors with the last property.

For any two vectors $x,y \in \R^k$, we say that {\it $y$ refines
$x$} (or that {\it $x$ is refined by $y$}), and write $x \preceq
y$, if $y_{i}=y_{j}$ implies $x_{i}=x_{j}$ for any $i,j \in \N_k$.

\begin{lemma}
\label{lem-3} Let $T_1$, $T_2$,...,$T_s$ be $s$-tensors on $\R^n$
with the following two properties:
\begin{enumerate}
\item for every multi index $(i_1,...,i_s)$ we have
$$
(T_1^{i_1...i_s}, T_2^{i_1...i_s},...,T_s^{i_1...i_s}) \preceq
(i_1,...,i_s).
$$
\item Whenever the multi index $(i_1,...,i_s)$ has pairwise
distinct entries, then $T_l^{i_1...i_s} = 0$ for all $l$.
\end{enumerate}
For any two multi indexes $(i_1,...,i_s)$, $(j_1,...,j_s)$ with
entries from the set $\N_n$ we define the $s \times s$ matrix,
$\Delta(\substack{i_1...i_s \\ j_1...j_s})$, as follows:
\[
\Delta(\substack{i_1...i_s \\ j_1...j_s})^{pq} = \left\{
\begin{array}{rl}
\delta_{i_pj_q}, & \mbox{ if } q < s, \\
T_{p}^{i_1...i_s} \delta_{i_pj_q}, & \mbox{ if } q = s,
\end{array}
\right.
\]
then $\det(\Delta(\substack{i_1...i_s \\ j_1...j_s})) = 0$.
\end{lemma}

\begin{proof}
Fix a multi index $(i_1,...,i_s)$. If it has pairwise distinct
entries, then by the second property, the last column in the
matrix $\Delta(\substack{i_1...i_s \\ j_1...j_s})$ will be zero.
If two of its entries $i_{s_1}$ and $i_{s_2}$ are equal, then by
the first property, row $s_1$ will be equal to row $s_2$ and again
the determinant will be zero. \hfill \qed
\end{proof}

\vspace{-0.4cm}

\begin{note} \rm The second condition in the preceding lemma in
trivially satisfied if $s > n$.
\end{note}

\begin{proposition}
\label{k-p} Suppose that $T_1$, $T_2$,...,$T_s$ be $s$-tensors on
$\R^n$ with the two properties in Lemma~\ref{lem-3}, then the
following identity holds:
\begin{equation}
\sum_{\sigma \in P^{s-1}} \sum_{l=1}^{s} \signum(\sigma_{(l)})
\Diag^{\sigma_{(l)}} T_{l} =0.
\end{equation}
\end{proposition}

\begin{proof}
Fix an arbitrary multi index $(\substack{i_1...i_{s}
\\ j_1...j_{s}})$. In the following calculation, the third
equality uses the facts $P^{s} = \{\sigma_{(l)} \, | \, \sigma \in
P^{s-1}, \,\, l \in \N_{s}\}$ and $\sigma_{(l)}(l)=s$ for $l \in
\N_{s}$, $\sigma \in P^{s-1}$.
\begin{align*}
\big( \sum_{\sigma \in P^{s-1}} \sum_{l=1}^{s}
\signum(\sigma_{(l)}) \Diag^{\sigma_{(l)}} T_{l} \big)^{\substack{i_1...i_{s} \\
j_1...j_{s}}} &= \sum_{\sigma \in P^{s-1}} \sum_{l=1}^{s}
\signum(\sigma_{(l)})T_{l}^{i_1...i_{s}}
\delta_{i_1j_{\sigma_{(l)}}(1)} \cdots
\delta_{i_{s}j_{\sigma_{(l)}}(s)} \\
&= \det(\Delta(\substack{i_1...i_{s}
\\ j_1...j_{s}})) \\
&= 0,
\end{align*}
where $\Delta(\substack{i_1...i_{s}
\\ j_1...j_{s}})$ is defined in Lemma~\ref{lem-3} for the
family of tensors $T_{1}, T_{2},...,T_{s}$. \hfill \qed
\end{proof}

Recall the definitions of the linear maps $T \in T^{k,n}
\rightarrow T^{(l)}_{\mbox{\rm \scriptsize out}} \in T^{k+1,n},
\mbox{ for } l \in \N_k$. For convenience in this section we also
define
$$
T^{(k+1)}_{\mbox{\rm \scriptsize out}} \equiv 0, \mbox{ for every
} T \in T^{k,n}.
$$

\begin{lemma}
Suppose that $T \in T^{k,n}$ is a symmetric tensor such that
$T^{i_1...i_k}=0$ for any multi index $(i_1,...,i_k)$ with
pairwise distinct entries. Then the $k+1$ $(k+1)$-tensors
$$
T^{(1)}_{\mbox{\rm \scriptsize out}}, T^{(2)}_{\mbox{\rm
\scriptsize out}},...,T^{(k)}_{\mbox{\rm \scriptsize out}},
T^{(k+1)}_{\mbox{\rm \scriptsize out}}
$$
satisfy the two conditions in Lemma~\ref{lem-3} (with $s=k+1$).
\end{lemma}

\begin{proof}
To show the first property, fix a multi index $(i_1,...,i_{k+1})$
and suppose that $i_{p}=i_{q}$ for some $1 \le p < q \le k+1$. If
$q=k+1$, then $\big(T^{(p)}_{\mbox{\rm \scriptsize out}}
\big)^{i_1...i_{k+1}} = \big(T^{(q)}_{\mbox{\rm \scriptsize
out}}\big)^{i_1...i_{k+1}} = 0$. If $q < k+1$, then
\begin{align*}
\big(T^{(p)}_{\mbox{\rm \scriptsize out}} \big)^{i_1...i_{k+1}} &=
\frac{T^{i_1...i_{p-1}i_{k+1}i_{p+1}...i_{q-1}i_{q}i_{q+1}...i_k}
- T^{i_1...i_k}}{\mu_{i_{k+1}} - \mu_{i_p}} \\
&=\frac{T^{i_1...i_{p-1}i_{q}i_{p+1}...i_{q-1}i_{k+1}i_{q+1}...i_k}
- T^{i_1...i_k}}{\mu_{i_{k+1}} - \mu_{i_p}} \\
&=\frac{T^{i_1...i_{p-1}i_{p}i_{p+1}...i_{q-1}i_{k+1}i_{q+1}...i_k}
-T^{i_1...i_k}}{\mu_{i_{k+1}}
- \mu_{i_q}} \\
&=\big(T^{(q)}_{\mbox{\rm \scriptsize out}} \big)^{i_1...i_{k+1}},
\end{align*}
where in the second equality we used the fact that $T$ is
symmetric, while in the third we used $i_{p}=i_{q}$.

The verification of the second condition in Lemma~\ref{lem-3}
follows immediately from the fact that $T$ has that property.
\hfill \qed
\end{proof}

Using the fact that $T^{(k+1)}_{\mbox{\rm \scriptsize out}} \equiv
0$, we obtain the next corollary from Proposition~\ref{k-p}.

\begin{corollary}
Suppose that $T \in T^{k,n}$ is a symmetric tensor such that
$T^{i_1...i_k}=0$ for any multi index $(i_1,...,i_k)$ with
pairwise distinct entries. Then we have the following identity:
\begin{equation}
\sum_{\sigma \in P^{k}} \sum_{l=1}^{k} \signum(\sigma_{(l)})
\Diag^{\sigma_{(l)}} T^{(l)}_{\mbox{\rm \scriptsize out}} =0.
\end{equation}
\end{corollary}

\section{Lifting of a tensor determined by the cycle type of a permutation}

The goal of this final section is to prove a result in the spirit
of Corollary~\ref{jan11a}.

It is well known that every permutation $\nu$ on $\N_k$ has a
unique decomposition into a product of disjoint cycles. Denote by
$s$ the number of disjoint cycles. These cycles partition the set
of integers, $\N_k$, in a natural way: two integers $j,i \in N_k$
are in the same partition if $\nu^l(i)=j$ for some $l$. We
enumerate the sets in the partition in the natural way: let
$I_{\nu, 1}$ be the set of the partition that contains the integer
$1$; let $I_{\nu, 2}$ be the set that contains the smallest
integer not in $I_{\nu, 1}$; let $I_{\nu, 3}$ be the set
containing the smallest integer not in $I_{\nu, 1} \cup I_{\nu,
2}$, and so on.

Take a vector $x \in \R^k$. We will say that the {\it permutation
$\nu$ refines vector $x$} (or that $x$ is refined by $\nu$), and
write $x \preceq \nu$, if $x_{l}=x_{\nu(l)}$ for every
$l=1,2,...,k$.  In order to know a vector $x$, refined by $\nu$,
it is enough to know the value of one coordinate with index from
every cycle of $\nu$. In other words, for a fixed $x$, refined by
$\nu$, the vector $(p_1,...,p_s)$ defined by
$$
p_l := x_i, \mbox{ for every } l \in \N_s \mbox{ and some }i \in
I_{\nu, l},
$$
completely {\it specifies $x$ given $\nu$}. The ordering on the
cycles of $\nu$, we agreed on above, makes sure that this is a
one-to-one correspondence between $\R^s$ and the set $\{x \in \R^k
\, | \, x \preceq \nu\}$.

Let now $T$ be a tensor from $T^{s,n}$. (Notice that the dimension
of the tensor is equal to the number of cycles of $\nu \in P^k$.)
Clearly $s \le k$ with equality if, and only if, $\nu$ is the
identity permutation. We define $T^{\nu}$ to be a tensor from
$T^{k,n}$ defined component wise by
$$
(T^{\nu})^{i_1...i_k} = \left\{
\begin{array}{ll}
T^{p_1...p_s}, & \mbox{ if } (i_1,...,i_k) \preceq \nu, \\
0, & \mbox{ otherwise,}
\end{array}
\right.
$$
where $(p_1,...,p_s)$ is the vector that specifies $(i_1,...,i_k)$
given $\nu$. Informally speaking, $T^{\nu}$ has the tensor $T$
placed on the diagonal ``subspace'' defined by the coordinate
equalities $\{i_l = i_{\nu(l)} \, | \, l \in \N_k \}$ and zeros
every where else.

When $\tau_l$ is the transposition $(l,k+1)$ on $\N_{k+1}$, this
definition coincides with the definition of $T^{\tau_l}$ given by
Equation~(\ref{perm-lift}).

Fix a vector $\mu \in \R^n$. For the rest of this section, when a
block-constant tensor or a block-diagonal matrix is mentioned, it
will be understood that the blocks are determined by the vector
$\mu$ as explained in Section~\ref{NAPR}.  Let $M$ be a symmetric
matrix and $h \in R^n$ be such that
\begin{equation}
\label{fedun} U^T M_{\mbox{\scriptsize \rm in}} U = \Diag h,
\end{equation}
for some block-diagonal orthogonal matrix $U \in O^n$.

\begin{lemma}
For any block-constant tensor, $T$, on $\R^n$ we have
$$
T[h] = T[\diag M],
$$
where vector $h$ and matrix $M$ are defined by
Equation~(\ref{fedun}).
\end{lemma}

\begin{proof}
Suppose that $T$ is a $k$-tensor. The proof is a direct
calculation from the definitions:
\begin{align*}
(T[h])^{i_1...i_{k-1}} &= \sum_{i_k=1}^n T^{i_1...i_{k-1}i_{k}}h^{i_k} \\
&= \sum_{l=1}^{r} T^{i_1...i_{k-1}\iota_l} \sum_{i \in I_{l}} h^{i} \\
&= \sum_{l=1}^{r} T^{i_1...i_{k-1}\iota_l} \sum_{i \in I_{l}} M^{ii} \\
&= \sum_{i_k=1}^n T^{i_1...i_{k-1}i_{k}}M^{i_ki_k} \\
&= (T[\diag M])^{i_1...i_{k-1}}. \\[-1.4cm]
\end{align*}
\hfill \qed
\end{proof}

%The following Proposition, complements Theorem~3.3 in \cite{}.

\begin{proposition}
\label{prop-1} Let $\nu$ be a permutation on $\N_{k+1}$ with $s$
disjoint cycles such that $\nu(k+1)=k+1$. Then for any
block-constant tensor $T$ in $T^{s,n}$ we have the identity:
\begin{equation}
\label{last-eqn} T^{\nu}[h] = T^{\nu}[\diag M],
\end{equation}
where vector $h$ and matrix $M$ are defined by
Equation~(\ref{fedun}).
\end{proposition}

\begin{proof}
If the multi index $(i_1,...,i_{k},i_{k+1})$ is refined by $\nu$,
then $(i_1,...,i_{k},j)$ is also refined by $\nu$ for any $j \in
\N_n$. Moreover, the vector that specifies it (given $\nu$) will
look like $(p_1,p_2,...,p_{s-1},i_{k+1})$ because $\nu(k+1)=k+1$
and the cycle containing the integer $k+1$ has one element.
Therefore, as $(i_1,...,i_{k},i_{k+1})$ goes over all possible
multi indexes of dimension $k+1$, refined by $\nu$,
$(p_1,p_2,...,p_{s-1},i_{k+1})$ will go over all multi indexes of
dimension $s$. This correspondence is one-to-one. Fix an arbitrary
multi index, $(i_1,...,i_{k},i_{k+1})$. If it is not refined by
$\nu$, then
$$
(T^{\nu}[h])^{i_1...i_{k}} = \sum_{i_{k+1}=1}^{n}
(T^{\nu})^{i_1...i_ki_{k+1}}h^{i_{k+1}} = 0,
$$
and similarly the right-hand side of (\ref{last-eqn}) is equal to
zero. If the multi index is refined by $\nu$, then using the
previous lemma we compute:
\begin{align*}
(T^{\nu}[h])^{i_1...i_k} &= \sum_{i_{k+1}=1}^{n}
(T^{\nu})^{i_1...i_{k}i_{k+1}}h^{i_{k+1}} \\
&= \sum_{i_{k+1}=1}^{n} T^{p_1...p_{s-1}i_{k+1}}h^{i_{k+1}} \\
&= (T[h])^{p_1...p_{s-1}} \\
&= (T[\diag M])^{p_1...p_{s-1}} \\
&= \sum_{i_{k+1}=1}^{n} T^{p_1...p_{s-1}i_{k+1}}(\diag M)^{i_{k+1}} \\
&= (T^{\nu}[\diag M])^{i_1...i_k}.
\end{align*}

\vspace{-0.9cm}

\hfill \qed
\end{proof}

The following result complements Corollary~\ref{jan11a}. It
follows by combining Proposition~\ref{prop-1} and
Corollary~\ref{jan11abc}.

\begin{corollary}
Let $\nu$ be a permutation on $\N_{k+1}$ with $s$ disjoint cycles
such that $\nu(k+1)=k+1$. Then for any block-constant tensor $T$
in $T^{s,n}$ and permutation $\sigma$ on $\N_k$ we have the
identity:
\begin{equation}
\Diag^{\sigma} \big(T^{\nu}[h] \big) = \big(
\Diag^{\sigma_{(k+1)}} T^{\nu} \big)[M],
\end{equation}
where vector $h$ and matrix $M$ are defined by
Equation~(\ref{fedun}).
\end{corollary}

%***********************************************************************
                           %END TEXT
%***********************************************************************

%\bibliographystyle{plain}
%\bibliography{master}
%        % Here is the other way to do a bibliography: it uses BibTeX.
%        % You will need to have set up your references in a
%        % special format in a file with extension .bib (eg, myrefs.bib)
%        % For details, see the LaTeX references listed in ACS "how-to"
%        % handout T-5, "Use TeX Resources" or type  man bibtex
%        % on Fraser, Monashee, Beaufort or Selkirk.
%
%\vfill\eject

%************************************************************************
\end{document}